\documentclass[10pt]{amsart}
\usepackage[english]{babel}
\usepackage{amssymb}
\usepackage{amsmath}
\newcommand{\D}{\mathcal{D}}
\newcommand{\I}{\mathcal{I}}
\newcommand{\J}{\mathcal{J}}
\newcommand{\E}{\mathcal{E}}
\newcommand{\F}{\mathcal{F}}
\newcommand{\G}{\mathcal{G}}
\newcommand{\h}{\mathcal{H}}
\newcommand{\s}{\mathcal{S}}
\newcommand{\A}{\mathcal{A}}
\newcommand{\M}{\mathcal{M}}
\newcommand{\N}{\mathcal{N}}
\newcommand{\Pe}{\mathcal{P}}
\newcommand{\B}{\mathcal{B}}
\newcommand{\comp}{\mathcal{C}}
\newcommand{\sts}{\mathcal{O}_X}

\newcommand{\lb}{\mathcal{L}}
\newcommand{\torus}{\mathcal{T}}

\newcommand{\can}{\overline{\phantom{x}}}

\newtheorem{lemma}{Lemma}
\newtheorem{theorem}{Theorem}
\newtheorem{proposition}{Proposition}
\newtheorem{corollary}{Corollary}

\theoremstyle{definition}

\newtheorem{example}{Example}

\newtheorem{remark}{Remark}



\title{Jordan algebras over algebraic varieties\\
}

\date{30.8.2007}
\author{S. Pumpl\"un}
\email{susanne.pumpluen@nottingham.ac.uk}
\address{School of Mathematics\\
University of Nottingham\\
University Park\\
Nottingham NG7 2RD\\
United Kingdom
}

\begin{document}

\begin{abstract}
We construct Jordan algebras over a locally ringed space $(X,\sts)$ using generalizations of the Tits process and the
first Tits construction by Achhammer [Ach1, 2]. Some general results on the structure of these algebras are obtained.
Examples of Albert algebras over a Brauer-Severi variety with associated central simple algebra of degree 3
are given.
\end{abstract}

\maketitle

\section*{Introduction}

Let $R$ be a unital commutative ring of scalars. The theory of generically algebraic Jordan algebras over a
unital commutative base ring $R$ which are finitely generated projective as $R$-modules was developed by Loos [L2].
 Generically algebraic Jordan algebra of degree 3 over $R$ were constructed from a cubic norm
structure [L2, 3.9] (see McCrimmon [M1, Theorem 1] for base fields). These algebras play an important role in Jordan
theory since every simple exceptional Jordan algebra
 over a field is of this type. By McCrimmon and Tits, any exceptional  algebra  over a field, also called an Albert algebra,
  can be constructed from a special
cubic norm structure, called the first or second Tits construction (cf. Jacobson [J1]).
In [P-S1, 2], Petersson-Racine  defined Tits structures $(B,N,A)$ over $R$, where $B$ is
an algebra with involution $*$ carrying a cubic norm structure $N$, $A\subset H(B,*)$ is an ample subspace and
several properties connecting $B$, $N$ and $A$ hold. If these conditions are satisfied, the cubic norm
structure on $A$ can be extended to one on $A\oplus B$ which thus becomes a cubic Jordan algebra.
 In the following this construction which generalizes the second Tits construction
 is called the {\it classical Tits process}.

In his PhD thesis [Ach1, 2], Achhammer developed a generalized Tits process for algebras defined
over a locally ringed space $(X,\sts)$
which, applied to algebras over arbitrary rings,  generalized  the classical Tits process and, in particular,
the classical first Tits construction. He obtained some general results on Albert algebras and gave
examples of Albert algebras over projective space. Independently, Albert algebras over integral schemes  were investigated by
Parimala-Suresh-Thakur [P-S-T1, 2].
In [P-S-T2],  generalized first and second Tits constructions for Albert algebras
over a domain $R$ such that $2,3\in R^\times$ were introduced.
Albert algebras over $\mathbb{A}_k^2$ were constructed which do neither arise from a generalized first nor second
 Tits construction [P-S-T2], proving that the situation over schemes differs from the one over fields.

We study Jordan algebras $\mathcal{A}$ over a locally ringed space $(X,\sts)$
using Achhammer's Tits process and  first Tits construction which are formulated more generally than the ones
given in [P-S-T1, 2].
In particular, we obtain Jordan algebras over $X$ by iterating the first Tits process.
As an example, we look at algebras over a Brauer-Severi variety associated with a central simple algebra
of degree 3. Jordan algebras over algebraic curves of genus zero and one are investigated in [Pu3].

 After establishing some basic notation
in Section 1, some general results on Jordan algebras and cubic forms with adjoint and base point are collected in Section 2. The Tits process and the first Tits construction
are introduced in Sections 3 and 4, formulated for algebras over locally ringed spaces. We end by giving examples of Albert algebras over Brauer-Severi varieties associated
with central simple algebras of degree 3 in Section 6.
Throughout the paper, we summarize and explain the relevant results from Achhammer's thesis [Ach1], which has not been published
(for parts of it, cf. [Ach2]). They will again be applied in [Pu3].

The standard terminology from algebraic geometry is used, see Hartshorne's book [H].
For the standard terminology on Jordan algebras, the reader is referred to the books by
McCrimmon [M2], Jacobson [J1] and Schafer [Sch].
 In the following, let $(X,\sts)$ be a locally ringed space.

\section{Preliminaries}

\subsection{Algebras over $X$}

 For $P \in X$ let $\mathcal{O}_{P,X}$ be the local ring of
$\mathcal{O}_X$ at $P$ and $m_P$ the maximal ideal of $\mathcal{O}_{P,X}$. The corresponding residue class
field is denoted by $k(P)=\mathcal{O}_{P,X}/m_P$. For an $\mathcal{O}_X$-module $\mathcal{F}$ the stalk of
$\mathcal{F}$ at $P$ is denoted by $\mathcal{F}_P$. $\mathcal{F}$ is said to have {\it full support} if
${\rm Supp}\,\mathcal{F}=X$; i.e., if $\mathcal{F}_P\not=0$ for all $P\in X$. We call
$\mathcal{F}$ {\it locally free of finite rank} if for each $P\in X$ there is an open neighborhood $U\subset X$
of $P$ such that $\mathcal{F}|_U=\mathcal{O}_U^r$ for some integer $r\geq 0$. The {\it rank} of
$\mathcal{F}$ is defined to be ${\rm sup}\{{\rm rank}_{\mathcal{O}_{P,X}}\mathcal{F}_P\,|\, P\in X\}$.
 The term ``$\mathcal{O}_X$-algebra" (or ``algebra over $X$'') always refers to nonassociative
$\mathcal{O}_X$-algebras which are unital and locally free of finite rank as $\mathcal{O}_X$-modules.
 An $\mathcal{O}_X$-algebra $\mathcal{A}$ is called {\it alternative} if $x^2y=x(xy)$ and
$yx^2=(yx)x$ for all sections $x,y$ of $ \mathcal{A}$ over the same open subset of $X$.
 An algebra $\mathcal{A}$ over $\mathcal{O}_X$ is called  {\it separable} if
 $\mathcal{A}(P)$ is a separable $k(P)$-algebra for all $P\in X$.

Recall that an associative  $\mathcal{O}_X$-algebra $\mathcal{A}$ is called an {\it Azumaya algebra} if
$\mathcal{A}_P\otimes_{\mathcal{O}_{P,X}} k(P)$ is a central simple algebra over $k(P)$ for all $P\in X$ [K].

\subsection{Composition algebras over $X$} (cf. [P1]) An $\mathcal{O}_X$-algebra $\mathcal{C}$ is called a {\it composition algebra} over $X$
if it has full support, and if there exists a quadratic form $N \colon \mathcal{C} \to  \mathcal{O}_X$ such that
the induced symmetric bilinear form $N(u,v)  =
N(u+v)-N(u)-N(v)$ is {\it nondegenerate}; i.e., it determines a module isomorphism
$\mathcal{C} \overset{\sim}{\longrightarrow} \mathcal{C}^\vee
=\mathcal{H}om (\mathcal{C},\mathcal{O}_X)$,
and such that $N(uv)=N(u)N(v)$ for all sections $u,v$ of $\mathcal{C}$ over the
same open subset of $X$.
 $N$ is uniquely determined by these conditions and called the {\it norm} of $\mathcal{C}$. It is denoted by $N_\mathcal{C}$.
Given an algebra $\mathcal{C}$ over $X$ and a quadratic form $N \colon \mathcal{C}
\to \mathcal{O}_X$,  $\mathcal{C}$ is a composition algebra over $X$ with
norm $N$ if and only if $\mathcal{C}_P$ is a composition
algebra over $\mathcal{O}_{P,X}$ with norm $N_P$ for all $P \in X$.
 Composition algebras over $X$ are invariant under base change and exist only in ranks 1, 2, 4 or 8.
A composition algebra of constant rank 2 (resp. 4 or 8) is called a {\it quadratic  \'etale algebra} (resp. {\it quaternion algebra} or
{\it octonion algebra}). A composition algebra over $X$ of constant rank is called {\it split}
 if it contains a composition subalgebra
isomorphic to $\mathcal{O}_X \oplus \mathcal{O}_X$.  There are several construction methods for composition algebras
over locally ringed spaces:  There  exists a Cayley-Dickson process ${\rm Cay}(\D,\Pe,N_\Pe)$
which is described in [P1]. If $2\in H^(X,\sts)$, every quaternion algebra $\mathcal{C}$ over $X$ can be built out of
a locally free $\sts$-module $\M$ of constant rank 3 with trivial determinant, carrying a nondegenerate
quadratic form $N_0$, cf. [Pu1, 2.7]. We write $\mathcal{C}={\rm Quat}(\M,N_0)$ and observe that
$\mathcal{C}=\sts\oplus \M$ as $\sts$-module.

\subsection{Forms of higher degree over $X$}

Let $d$ be a positive integer.  Let $\mathcal{M}$, $\N$ be $\mathcal{O}_X$-modules which are locally free of finite rank.
 In general, that is if there are no restrictions on the global sections $H^0(X,\sts)$,
 polynomial maps from $\M$ to $\N$ can be defined
analogously as it was done in [R], [L1, \S 18] or [L2] for polynomial maps between modules over rings, see  [Ach1].
In that case forms of degree $d$ are usually identified with the induced polynomial maps.

If $d!\in H^0(X,\mathcal{O}_X^\times)$, which will be the case considered in large parts of the paper,
 this definition is equivalent to the following:
A {\it form of degree $d$} $N:\mathcal{M}\to \sts$  over $X$ is a map
 such that $N(a x)=a^d N(x)$ for all sections $a$ in $ \sts$, $x$ in $\M$, where the map
 $$\theta : \M \times\dots \times \M \to \sts$$
 defined by
 $$\theta(x_1,\dots,x_d)=\frac{1}{d!}\sum_{1\leq i_1<\dots<i_l\leq d}(-1)^{d-l}N(x_{i_1}+\dots+x_{i_l})$$
 is a $d$-linear form over $\sts$ (the range of summation of $l$ being $1\leq l\leq d$).
  $\theta$ is called the {\it  symmetric $d$-linear form}
associated  with $N$ and $(M,\theta)$ a {\it $d$-linear space}.
 Obviously, $N(x)=\theta(x,\dots,x)$ for all sections $x$ of $\mathcal{M}$ over the same open subset of $X$. We
   do not distinguish between a form of degree $d$ and its associated symmetric $d$-linear form $\theta$.

Let
$$ \M\to {\rm \h om}_{X}(\M\otimes\dots\otimes \M,\sts)$$
 (($d-1$)-copies of $\M$) be  the map
$$\theta_{x_1}(x_2\otimes\dots\otimes x_d)=\theta(x_1,x_2,\dots,x_d).$$
$N:\mathcal{M}\to \sts$ (or its associated $d$-linear form $\theta$) is called {\it  nondegenerate}, if
 $\theta_{x_1}\otimes k(P)$ is injective for all $P\in X$.
Thus $N$ is nondegenerate if and only if $N(P):\mathcal{M}(P)\to k(P)$ is nondegenerate in the classical sense for
 all $P\in X$. This concept of nondegeneracy is invariant under base change.
 In the present paper, only nondegenerate forms will be investigated, unless specified otherwise.

Two $d$-linear spaces $(\M_i,\theta_i)$, $i=1,2$ are called {\it isomorphic} (written
$(\M_1,\theta_1)\cong (\M_2,\theta_2)$ or just $\theta_1\cong\theta_2$)
if there exists an $\sts$-module isomorphism $f:\M_1\to \M_2$ such that
$\theta_2(f(x_1),\dots,f(x_d))=\theta_1(x_1,\dots,x_d)$  for all sections $x_1,\dots,x_d$ of
$\mathcal{M}_1$ over the same open subset of $X$.

The {\it orthogonal sum} $(\M_1,\theta_1)\perp (\M_2,\theta_2)$ of $(\M_i,\theta_i)$, $i=1,2$, is defined to be the $\sts$-module
 $\M_1\oplus \M_2$ together with the $d$-linear form
$(\theta_1 \perp \theta_2)(u_1+x_1,\dots,u_d+x_d)=\theta_1(u_1,\dots,u_d)+\theta_2(x_1,\dots,x_d)$.
A $d$-linear space $(\M,\theta)$ is
called {\it decomposable}, if $(\M,\theta)\cong (\M_1,\theta_1)\perp (\M_2,\theta_2)$
for two non-zero $d$-linear spaces $(\M_i,\theta_i)$, $i=1,2$.
A non-zero $d$-linear space $(\M,\theta)$ is
called {\it indecomposable} if it is not decomposable. We  distinguish between indecomposable ones and
{\it absolutely} indecomposable ones;
i.e., $d$-linear spaces which stay indecomposable under base change. Forms of degree 3 are also called {\it cubic forms}.

\begin{lemma} Suppose $3\in H^0(X,\sts)$. Every line bundle $\mathcal{N} \in \, _3 {\rm Pic} X$ carries a
nondegenerate cubic form which is uniquely determined up to multiplication with an invertible element in $H^0(X,\sts)$.
\end{lemma}

\begin{proof} If $\beta:\N\otimes\N\otimes\N\to\sts$ is an isomorphism then
$$N(w)=\beta(w\otimes w\otimes w)$$
is a nondegenerate cubic form on $\N$. $N$ is a norm on $\N$ and hence is uniquely determined
up to multiplication with an invertible element in $H^0(X,\sts)$ [Ach1, 2.11].
\end{proof}

\section{Jordan algebras and cubic forms with adjoint and base point over $X$}

\subsection{Jordan algebras over $X$}
Let $\J$ be a locally free $\sts$-module of
finite rank.  $( \mathcal{J},U,1)$ with $1\in H^0(X,\mathcal{J})$ is a {\it Jordan algebra} over $X$ if:
\begin{enumerate}
\item The $U$-operator $U:\mathcal{J}\to\mathcal{E}nd_{\mathcal{O}_X}(\mathcal{J}), x\to U_x$ is a quadratic map;
\item $U_1=id_\mathcal{J}$;
\item $U_{U_x(y)}=U_x\circ U_y\circ U_x$ for all sections $x,y$ in $ \mathcal{J}$;
\item $U_{x}\circ U_{y,z}(x)=U_{x,U_x(z)}(y)$ for all sections $x,y,z$ in $ \mathcal{J}$;
\item for every commutative associative $\sts$-algebra $\sts'$, $\mathcal{J}\otimes \sts'$ satisfies $(3)$ and $(4)$.
\end{enumerate}
In general, we  write  $\J$ instead of $( \mathcal{J},U,1)$ [Ach1].

 An  $\mathcal{O}_X$-algebra $\mathcal{J}$ is called an {\it Albert algebra}
 if $\mathcal{J}(P)=J_P\otimes k(P)$ is an Albert algebra over $k(P)$
for all $P\in X$.

\begin{remark}  (i)  With the usual definition of quadratic Jordan algebras over rings [J2, 1.3.4] we have:
Let  $ \mathcal{J}$ be a locally free $\mathcal{O}_X$-module
of finite rank, $1\in H^0(X,\mathcal{J})$ and $U:\mathcal{J}\to\mathcal{E}nd_{\mathcal{O}_X}(\mathcal{J}),
 x\to U_x$ a quadratic map. Then $( \mathcal{J},U,1)$
 is a Jordan algebra over $X$, iff $( \mathcal{J}(V),U(V),1|_V)$ is a Jordan algebra over $\sts(V)$
 for all open subsets $V\subset X$, iff $( \mathcal{J}_P,U_P,1_P)$ is a Jordan algebra over ${\sts}_P$ for all $P\in
X$.\\
(ii) Let $\A$ be a unital associative algebra over $X$. Then $\A^+=(\A,U,1)$ with $U:\A\to
\mathcal{E}nd_{\mathcal{O}_X}(\mathcal{A}), x\to U_x(y)=xyx$, $1=1_\A$, is a Jordan algebra over $X$.
\\(iii)  There is a canonical equivalence between the category of Jordan algebras
over $R$, which are finitely generated projective as $R$-modules, and the category of Jordan
algebras over the
affine scheme $Z = {\rm Spec} \, R$, which are locally free as $\sts$-modules,
given by the functor $(J,U,1) \longrightarrow (\widetilde{J}, \widetilde{U}, \widetilde{1})$.\\
(iii) Jordan algebras are invariant under base change: If $\sigma:X'\to X$ is a morphism of locally ringed spaces and
$( \mathcal{J},U,1)$ a Jordan algebra over $X$ then $\sigma^*( \mathcal{J},U,1)=(\sigma^* \mathcal{J},\sigma^*U,1)$
is a Jordan algebra over $X'$. Moreover, if $( \mathcal{J}',U',1')$ a Jordan algebra over $X'$ then
$\sigma_*( \mathcal{J}',U',1')=(\sigma_* \mathcal{J}',\sigma_*U',1')$
is a Jordan algebra over $X$.
\end{remark}

 If $\mathcal{J}$ is a Jordan algebra over a scheme $(X,\mathcal{O}_X)$,
then $\mathcal{J}$ is an Albert algebra if and only if there is a covering
$V_i\to X$ in the flat topology on $X$ such that $\mathcal{J}\otimes \mathcal{O}_{V_i}\cong \mathcal{H}_3({\rm Zor}( \mathcal{O}_{V_i}))$
where ${\rm Zor}(  \mathcal{O}_{V_i})$ denotes the octonion algebra of Zorn vector matrices over  $\mathcal{O}_{V_i}$
and  $\mathcal{H}_3({\rm Zor}( \mathcal{O}_{V_i}))$ is the Jordan algebra of 3-by-3 hermitian matrices
with entries in ${\rm Zor}( \mathcal{O}_{V_i})$ and scalars $ \mathcal{O}_{V_i}$  on the diagonal [Ach1, 1.10].
This terminology is compatible with the one used in [P-S-T1], as it was already observed in [P2, Section 2].

\begin{lemma} Let $X$ be a scheme over the affine scheme $Y={\rm Spec}\,R$ and suppose $H^0(X,\mathcal{O}_X)=R$.
Then a Jordan algebra (resp., an Azumaya algebra) $\mathcal{J}$ over $X$ is defined over $R$ providing it is globally free as an
$\mathcal{O}_X$-module.
\end{lemma}

The proof is analogous to [P1, 1.10].

\begin{remark} (analogous to [P1, 4.6]) Let $k$ be a field and $Y$ any $k$-scheme such that $H^0(Y,\mathcal{O}_Y)=k$.\\
(i) Let $A$ be
an Azumaya algebra over $k$. Write $\A$ for the base change from $A$ to $Y$. The canonical map $\rho:A\to H^0(Y,\A)$
is an isomorphism of $k$-algebras by [K, p.~143], since both algebras have the same dimension and $A$ is central simple.
 In particular, we recover $A$ from $\A$ by passing to global sections.
\\ (ii) Let $J$ be a simple Jordan algebra over $k$. Write $\J$ for the base change from $J$ to $Y$.
The canonical map $\rho:J\to H^0(X,\J)$ is an isomorphism of $k$-algebras: both algebras have the same dimension,
$J$ is a Jordan algebra and so is $H^0(Y,\A)$. Since $J$ is simple, the kernel of $\rho$ must be zero and
  $\rho$ is injective.
\end{remark}

\begin{proposition} Let $R$ be a ring and $Y$ an $R$-scheme such that $H^0(Y,\mathcal{O}_Y)=R$. Write $\sigma:Y\to Z={\rm Spec}\,R$
for the structure morphism of $Y$.\\
(i) Let $\A$ be an Azumaya algebra over $Y$
and $A_0\subset H^0(Y,\A)$ an Azumaya subalgebra over $R$.
Then the $\mathcal{O}_Y$-submodule $\A_0$ of $\A$ generated by the global sections in $A_0$ is an Azumaya subalgebra
of $\A$ defined over $R$. More precisely, there is a natural isomorphism from $\A_0$ to the base change of $A_0$
from $R$ to $Y$ (i.e., $\A_0\cong A_0\otimes \mathcal{O}_Y$).
\\(ii)
Let $\J$ be a Jordan algebra over $Y$  and $J_0\subset H^0(Y,\J)$ a Jordan subalgebra
over $R$ such that $\sigma^*(\widetilde{J_0})_P\otimes k(P)$ is simple for every $P\in Y$.
Then the $\mathcal{O}_Y$-submodule $\J_0$ of $\J$ generated by the global sections in $J_0$ is a Jordan subalgebra
of $\J$ defined over $R$. More precisely, there is a natural isomorphism from $\J_0$ to the base change of $J_0$
from $R$ to $Y$.
\end{proposition}

The proof is analogous to [P1, 5.2].

\subsection{Cubic forms with adjoint and base point}  Let $\J$ be an $\sts$-module.
A tripel $(N,\sharp,1)$ is a cubic form with adjoint and base point on $\J$ if $N:\J\to \sts$ is a cubic form,
$\sharp:\J\to\J$ a quadratic map and $1\in H^0(X,\J)$, such that
$$\begin{array}{l}
 x^{\sharp\,\sharp}=N(x)x,\\
 T(x^\sharp,y)=D_yN(x) \text{ for } T(x,y)=-D_xD_y {\rm log} N(1),\\
 N(1)=1,\, 1^\sharp=1,\\
 1\times y=T(y)1-y \text{ with } T(y)=T(y,1),\, x\times y=(x+y)^\sharp-x^\sharp-y^\sharp
\end{array}$$
for all sections $x,y$ in $\J$ over the same open subset of $X$ ([P-R1], [Ach1]). Here, $D_yN(x)$ denotes the
directional derivative of $N$ in the direction $y$, evaluated at $x$.
The symmetric bilinear form $T$ is called the {\it trace form} of $\J$.

Every cubic form with adjoint and base point $(N,\sharp,1)$ on a locally free $\sts$-module $\J$ of finite rank
defines a unital Jordan algebra structure $\J(N,\sharp,1)=(\J,U,1)$ on $\J$ via
$$U_x(y)=T(x,y)x-x^\sharp\times y$$
for all sections $x,y$ in $\J$, where the identities given in
[P-R1, p.~213] hold for all sections in $\J$. A section $x\in\J$ is invertible iff $N(x)$
is invertible. In that case, $x^{-1}=N(x)^{-1}x^\sharp$ ([P-R1], [Ach1]).

\begin{theorem} ([Ach1, 1.13, 1.14]) Let $\mathcal{J}$ be a Jordan algebra over $X$. Let $N$, $S$ and $T$ be
a  cubic, quadratic and linear form from $\mathcal{A}$ to $\sts$ such that
$$x^3-T(x)x^2+S(x)x^2-N(x)1=0$$
for all sections $x$ in $\J$. Suppose that for each $P\in X$ there exists an element $u\in \mathcal{J}(P)$ such that $1,u,u^2$ are linearly
independent over $k(P)$.  Let $M$, $Q$ and $L$ be some other cubic,  quadratic and linear forms from
$\mathcal{A}$ to $\mathcal{O}_X$
satisfying $$x^3-L(x)x^{2}+Q(x)x^{}-M(x)1=0$$ for all sections $x$ in $\mathcal{J}$.
Then $M=N$, $S=T_2$ and $L=T$.
\end{theorem}

 Over a field $k$ with infinitely many elements, such $1,u,u^2$ exist for every
 nondegenerate generically algebraic Jordan algebra of degree 3 over $k$. Furthermore,
  such elements exist for each Albert algebra and for every $J=A^+$, $A$ an Azumaya algebra of rank 9,
  over any field [Ach1, 1.13].

Here are some examples of cubic forms with adjoints translated from the well-known  setting of
cubic norm structures over rings (cf. [P-R1, p.~214 ff.]).

\begin{example} (i) Let $\A=\sts\times\sts\times\sts$, $N_\A(a_1,a_2,a_3)=a_1a_2a_3$,
$(a_1,a_2,a_3)^\sharp_\A=(a_2a_3,a_3a_1,a_1a_2)$ and $1=(1,1,1)$. Then $(N_\A,\sharp_\A,1)$ is a cubic
form with adjoint and base point on $\A$ and $\J(N_\A,\sharp_\A,1)=\A^+$.\\
(ii) Let $Q':\M\to \sts$ be a quadratic form on an $\sts$-module $\M$ with base point $1'$; i.e., $1'\in H^0(X,\M)$
satisfies $Q'(1')=1\in H^0(X,\sts)$. Put $\J=\sts\oplus \M$, $N(a,v)=aQ'(v)$, $(a,v)^\sharp=(Q'(v),a\bar{v})$,
where $\bar{v}=T'(v)1'-v$, $T'(v)=Q(1',v)$, and $1=(1,1')$. Let $T'(u,v)=Q'(v,\bar{u})$. Then
$\J(N,\sharp,1)=\sts^+\oplus \J(Q',1')$. Note that $\J(Q',1')$ has $wU_v=Q(v,\bar{w})v-Q'(v)\bar{w}$.

In particular, let $\comp'$ be a unital quadratic alternative algebra over $X$ with a scalar involution
$\bar{x}$ and with unit $1'$, norm $n'(x)1'=x\bar{x}$ and trace $t'(x)=x\bar{x}$
 (e.g., a composition algebra $\comp'$ with canonical involution $\can$).
Then $\B=\sts\oplus\comp$ is a unital alternative algebra over $X$ and $\B^+=\J(N,\sharp,1)$ with $(N,\sharp,1)$
defined using $n'$, $t'$, $\bar{\,}\,'$ as above. The resulting $(N,\sharp,1)$ is $\B$-admissible
with respect to the involution $(a,u)^*=(a,\bar{u}\,')$ (see 3.1 for the definition of $\B$-admissible). \\
(iii) Let $N,\sharp:\sts\to\sts$, $N(x)=x^3$, $a^\sharp=a^2$ and $1\in H^0(X,\sts)$. Then $\J(N,\sharp,1)=\sts^+$ and
 the associated trace form is $T(a,b)=3ab$.
\\ (iv) Let $\mathcal{C}$ be a unital quadratic alternative $\sts$-algebra with scalar involution $\bar{\,}\,'$
of constant rank over $X$.
Let $\Gamma={\rm diag}\, (\gamma_1, \gamma_2, \gamma_3)$ be a diagonal matrix with
$\gamma_i\in H^0(X,\mathcal{O}^\times_X)$. Let $\sigma:{\rm Mat}_3(\comp)\to {\rm Mat}_3(\comp)$
be the involution given by $\sigma(X)=\Gamma^{-1}\overline{ X}^t \Gamma$, bar denoting the entry-wise action of the canonical
involution of $\comp$.   Let $\mathcal{J}=\mathcal{H}_3(\comp,\Gamma)$ be
 the $\sts$-submodule of elements of the form
$$x=\sum_{i=1}^3 a_ie_i+\sum_{ijk}u_{i[jk]}$$
with $a_i\in \sts$, $e_i=e_{ii}$, $u_i\in\comp$, $u_{i[jk]}=u_ie_{jk}+\gamma_k^{-1}\gamma_j\bar{a_i}e_{kj}$
 where $\sum_{ijk}$ denotes the sum over all cyclic permutations of $(123)$.
Then
$$\begin{array}{l}
N(x)=a_1a_2a_3-\sum_{ijk}a_i\gamma_k^{-1}\gamma_jn_\comp(u_i)+t_\comp(u_1u_2u_3),\\
x^\sharp=\sum_{ijk}(a_ja_k-\gamma_k^{-1}\gamma_jn_\comp(u_i))+\sum_{ijk}(\gamma_j^{-1}\gamma_k
\overline{u_ju_k}-a_iu_i)_{[jk]},\\
1=\sum_{i=1}^3 e_i
\end{array}$$
is a cubic form  on $\J$ with adjoint and base point.
(If $2\in H^0(X,\sts)$, then $\mathcal{H}_3(\comp,\Gamma)$ is the sheaf of symmetric elements for $\sigma$.)
If $\Gamma= {\rm diag}\, (1, 1,1)$ we write $\mathcal{H}_3(\comp)=\mathcal{H}_3(\comp,\Gamma)$.
\\(v) Let $X$ be a scheme and $\J$ an Albert algebra over $X$. Then there is a cubic form $N:\J\to\sts$ with
adjoint $\sharp:\J\to\J$ and base point 1 such that $\J=\J(N,\sharp,1)$ [Ach1, 1.17].
\end{example}

\begin{lemma} Let $2\in H^0(X,\sts)$ be invertible. Let $\J_1=\h_3(\comp_1,\Gamma_1)$,  $\J_2=\h_3(\comp_2,\Gamma_2)$
be two Jordan algebras over $X$ with trace forms $T_1$, respectively $T_2$, and let $\comp_i\not=0$.\\
(i) Let $\comp_1$ be an octonion algebra or let $|k(P)|$ be infinite
for all $P\in X$ and $\comp_1(P)\cong \comp_2(P)$ for all $P\in X$. If $\J_1\cong\J_2$ then $T_1\cong T_2$.\\
(ii) If $T_1\cong T_2$ then $\J_1(P)\cong\J_2(P)$ and $\comp_1(P)\cong \comp_2(P)$ for all $P\in X$.
Let in particular $X$ be a proper scheme
over a perfect field. If the norm restricted to the subspace of trace zero elements in $\comp_1$ is indecomposable, then
$\comp_1\cong\comp_2$.
\end{lemma}

This follows from [KMRT, IX, (37.13)] applied to the residue class algebras. For examples of composition algebras, where
the norm restricted to the subspace of trace zero elements is indecomposable, see for instance
the octonion algebras over $\mathbb{A}_K^2$ constructed in [T] or
 the quaternion algebras over curves of genus one constructed in [Pu1].

\begin{example} Let $2\in H^0(X,\sts).$\\
 (i) $\J=\mathcal{H}_3(\mathcal{O}_X,\Gamma)=\{A\in {\rm Mat}_3(\mathcal{O}_X)\,|\,\sigma(A)=A  \}$ is a free
 $\mathcal{O}_X$-module of rank 6. If $X$ is a scheme over $Y={\rm Spec}\,R$
  and $H^0(X,\mathcal{O}_X)=R$, then $\mathcal{J}$ is defined over $R$ by Lemma 2.\\
(ii) If $\torus={\rm Cay}(\sts,\lb,N)$ is a torus over $X$,
 then $\mathcal{H}_3(\torus,\Gamma)$ has rank 9 and its underlying $\sts$-module structure is
 given by $\mathcal{H}_3(\torus,\Gamma)=\mathcal{O}^6_X\oplus  \lb^3$.
 For a split quadratic  \'etale algebra,
 $\mathcal{H}_3(\mathcal{O}_X\times \mathcal{O}_X)\cong{\rm Mat}_3(\mathcal{O}_X)^+$.
\end{example}

\begin{remark} Let $2\in H^0(X,\sts)$ be invertible. Let $\J$ be a Jordan algebra of constant rank greater
 than 3 containing a split cubic  \'etale algebra $\mathcal{D}$. Identify
$\mathcal{D}=\sts\times\sts\times\sts$. Since $\mathcal{D}$ is a cubic alternative algebra, so is
$D=H^0(X,\mathcal{D})$.
$D^+\subset H^0(X,\mathcal{J})$ has a complete orthogonal system of primitive idempotents $(c_1,c_2,c_3)$ which
  gives rise to a refined Peirce decomposition (similar as in the proof of [P1, 3.5])
  $$\J=\J_{1,1}\oplus \J_{1,2}\oplus \J_{1,3}\oplus \J_{2,1}\oplus \J_{2,2}\oplus \J_{2,3}\oplus
   \oplus \J_{3,1}\oplus \J_{3,2}\oplus \J_{3,3}$$
   of $\J$ with $\J_{i,i}=\sts c_i$ ($1\leq i\leq 3$) and
   $$\J_{i,j}=\J_{j,i}=\{x\in \J\,|\,xc_i=\frac{1}{2}x=xc_j \}$$
    ($1\leq i<j\leq 3$) [Sch, p.~100]. If $\J$ has rank 9, $\J_{i,j}$ is a line bundle for all $i,j$.
\end{remark}

\section{The Tits process}

\subsection{} We introduce the generalized Tits process for locally ringed spaces due to [Ach1]:
Let $(X,\sts')$ be a locally ringed space, $*:\sts'\to \sts'$ an involution, $\B$ an associative $\sts'$-algebra
and $*_\B$ an involution on $\B$ such that $*_\B|_{\sts'}=*$. To simplify the notation, we
write $*_\B$ also for the involution $*$ on $\sts'$. Let $(N_\B,\sharp_\B, 1)$ be a
cubic form with adjoint and base point on $\B$.
Assume that $(N_\B,\sharp_\B, 1)$ is $\B$-{\it admissible}.
That means,
\begin{enumerate}
\item $\B^+=J(N_\B,\sharp_\B,1)$ with $1\in H^0(X,\B)$ the unit element in $\B$, and $xyx=T_\B(x,y)x-x^{\sharp_\B}\times_\B
y$ for $x,y$ in $\B$;
\item $N_\B(xy)=N_\B(x)N_\B(y)$ for all $x,y$ in $\B$;
\item $N_\B(x^{*_\B})=N_\B(x)^{*_\B}$ for $x$ in $\B$
\end{enumerate}
[Ach1, 2.6]. Write $\B^\times$ for the sheaf of units
of $\B$. The pointed set of isomorphism classes ${\rm Pic}_{l}\B$ of locally free left $\B$-modules of rank one can be matched
canonically with $\check{H}^1(X,\B^\times)$ in the sense of noncommutative \u{C}ech cohomology [Mi]. The morphisms of group
sheaves $N_\B:\B^\times\to \sts'^\times$ and
$\sharp_\B:\B^\times\to (\B^{\rm op})^\times$ determine morphisms of pointed sets
$$N_\B: \check{H}^1(X,\B^\times)\to \check{H}^1(X,{\sts'}^\times)={\rm Pic}\,\sts'$$
and
$$\sharp_\B: \check{H}^1(X,\B^\times)\to \check{H}^1(X,(\B^{op})^\times).$$
After suitable identifications, we get a morphism
$$\sharp_\B: {\rm Pic}_{l}\B\to  {\rm Pic}_{l}\B^{\rm op}.$$
A cubic map $f:\Pe \to \E$ in the category of $\sts'$-modules is called {\it multiplicative} if
$$f(bw)=N_\B(b)f(w)$$
 for all sections $b$ in $\B$, $w$ in $\E$. Let $\F$ be a right $\B$-module.
 A quadratic map $g:\Pe \to \F$ in the category of $\sts'$-modules is called {\it multiplicative} if
$$g(bv)=g(v)b^{\sharp_\B}$$
 for all sections $b$ in $\B$, $v$ in $\Pe$. A multiplicative cubic map $N:\Pe\to N_\B(\Pe)$ is called a {\it norm}
 on $\Pe$ if $N$ is universal in the category of multiplicative cubic maps on $\Pe$.
 A multiplicative quadratic map $\sharp:\Pe\to \Pe^{\sharp_\B}$ is called an {\it adjoint} on $\Pe$ if $\sharp$
 is universal in the category of multiplicative quadratic maps on $\Pe$.

 By the usual gluing process (cf. [P1, 2.4]), it can be shown that cubic norms on $\Pe$ always exist and are unique
 up to an invertible factor in $H^0(X,\sts')$ [Ach1, 2.11]. Analogously,
 one can prove the following [Ach1, 2.13]:

Let $\Pe^\vee$ denote the right $\B$-module $\mathcal{H}om_\B(\Pe,\B)$.
Let $\Pe\in {\rm Pic}_{l}\B$ such that $N_\B(\Pe)\cong \sts'$ and let $N:\Pe\to \sts'$ be a cubic norm on $\Pe$.
Then
$$\Pe^{\sharp_B}\cong \Pe^\vee,\quad {\Pe^\vee}^{\sharp_B}\cong \Pe$$ and
$N_\B(\Pe^\vee)\cong\sts'$.
There exists a uniquely determined cubic norm $\check{N}:\Pe^\vee\to \sts'$ and
uniquely determined adjoints $\sharp:\Pe\to \Pe^\vee$ and $\check{\sharp}:\Pe^\vee\to \Pe$ such that
\begin{enumerate}
\item $\langle w,w^\sharp\rangle=N(w)1$;
\item $\langle \check{w}^{\check{\sharp}},\check{w}\rangle=\check{N}(\check{w})1$;
\item $w^{\sharp \, \check{\sharp}}=N(w)w$
\end{enumerate}
for all $w$ in $\Pe$, $\check{w}$ in $\Pe^\vee$. Moreover,
\begin{enumerate}
\item $\check{w}^{\check{\sharp}\,\sharp }=\check{N}(\check{w})\check{w}$;
\item $\langle w, \check{w}\rangle^{\sharp_\B}= \langle \check{w}^{\check{\sharp}}, w^{\sharp}\rangle$;
\item $N_\B(\langle w,w^\sharp\rangle)=N(w)\check{N}(\check{w})$;
\item $D_{w'}N(w)=T\B(\langle w', w^\sharp\rangle)$;
\item $D_{\check{w}'}\check{N}(\check{w})=T\B(\langle\check{w}^{\check{\sharp}},\check{ w}'\rangle)$;
\item $\langle w,\check{w}\rangle w= T_\B(\langle w,\check{w}\rangle)w-w^\sharp\check{\times}\check{w}$
\end{enumerate}
for all $w,w'$ in $\Pe$, $\check{w}$ in $\Pe^\vee$. The morphism of group sheaves
$$*_\B:\B^\times\rightarrow (\B^{\rm op})^\times$$
determines  morphisms
$$\begin{array}{l}
*_\B:{\rm Pic}_{l}\,\B=\check{H}^1(X,\B^\times)\rightarrow \check{H}^1(X,(\B^{op})^\times)={\rm Pic}_{l}\,\B^{\rm op},\\
*_\B:{\rm Pic}_{l}\,\B^{\rm op}=\check{H}^1(X,(\B^{\rm op})^\times)\rightarrow \check{H}^1(X,\B^\times)={\rm Pic}_{l}\,\B
\end{array}$$
of pointed sets.

\subsection{}
For a locally free left $\B$-module $\Pe$ (respectively, right $\B$-module) of rank one, let $\overline{\Pe}$ denote the opposite
module of $\Pe$  with respect to the involution $*_\B$ [K, I.(2.1)]. An isomorphism of left
(respectively, right)  $\B$-modules $j:\Pe \to \overline{{\Pe^*}^{\B}}$ is called an {\it involution } on $\Pe$.
Furthermore, we canonically identify the left $\B$-homomorphisms from $\Pe$ to $\overline{\Pe^\vee}$ with the
sesquilinear forms on $\Pe$ [K, I.(2.2)].

 A pair $(\A,\sts)$ consisting of a subsheaf of rings $\sts$
of $\sts'$ and an $\sts$-submodule $\A$ of $\B$ is called $\B$-{\it ample} if
\begin{enumerate}
\item $\sts\subset \mathcal{H}(\sts',*_{\B})$,
\item $rr^{*_\B}\in\sts$ for $r$ in $\sts'$,
\item $\A\subset \mathcal{H}(\B,*_{\B})$,
\item $1\in H^0(X,\A)$,
\item $bab^*\in \A$ for $a\in\A$, $b\in\B$,
\item $N_\B(\A)\subset\sts$; i.e., $N_\B|_\A:\A\to\sts$ is a cubic form over $\sts$,
\item $\A^{\sharp_\B}\subset \A$; i.e., ${\sharp_\B}|_\A: \A\to\A$ is a quadratic map over $\sts$.
\end{enumerate}
Let $(\A,\sts)$ be $\B$-ample and $\Pe$ be a locally free left $\B$-module of rank 1 with $N_\B(\Pe)\cong \sts'$.
\\ (i) If $\Pe^{*_\B}\cong\Pe^\vee$ and $N_\B(\Pe)\cong\sts'$, then a pair $(N, *)$ with $N:\Pe\to\sts'$ a norm on $\Pe$
and an involution $*:\Pe\to \overline{\Pe^\vee}$ on $\Pe$ is called $\A$-{\it admissible} if
\begin{enumerate}
\item $\langle w,w^*\rangle\in\A$,
\item $N_\B(\langle w,w^*\rangle)=N(w)N(w)^{*_\B}$
\end{enumerate}
for $w\in\Pe$.\\
(ii) $\Pe$ is called $\A$-{\it admissible} if there is a cubic norm $N:\Pe\to \mathcal{O}_X'$ and a nondegenerate
 $*_\B$-sesquilinear
form $h:\Pe\times\Pe\to \B$ (i.e., $h(aw,bv)=ah(w,v)b^{*_\B}$ and $h$ induces an $\B$-module isomorphism
$j_h:\Pe\to \overline{\Pe^\vee}$) such that
\begin{enumerate}
\item $h(w,w)\in\A$,
\item $N_\B(h(w,w))=N(w)N(w)^{*_\B}$
\end{enumerate}
for $w\in\Pe$.
Note that $\Pe^\vee\cong\Pe^{*_\B}$   and that therefore $j_h$ (denoted $*$ from now on) is an involution  on $\Pe$
such that
\begin{enumerate}
\item $\langle w,w^*\rangle\in\A$,
\item $N_\B(\langle w,w^*\rangle)=N(w)N(w)^{*_\B}$
\end{enumerate}
for $w\in\Pe$.

\begin{remark} (i) $(\mathcal{H} (\B,*_\B),\mathcal{H} (\sts',*_\B))$ is $\B$-ample. If $2\in H^0(X, \sts')$ is invertible or if there is an element
$\lambda\in H^0(X, \sts')$ such that $1=\lambda+\lambda^{*_\B}$, this is the only $\B$-ample pair ([P-R1], [Ach1, 2.17]).
\\ (ii) The definition of $\B$-admissible in [P-R1] also includes nonassociative, e.g. alternative,
algebras and thus is more general. However, the Tits process introduced by Achhammer needs $\B$ to be
associative in order to be able to use non-trivial left $\B$-modules $\Pe$ in the construction. There would be no sensible
notion of a left $\B$-module otherwise.
\end{remark}

\subsection{} Keeping the notations of 4.1, 4.2, let $\Pe$ be $\A$-admissible. Then $\check{*}=\overline{*^{-1}}:\Pe^\vee\to \overline{\Pe}$ is an involution
on $\Pe$ and
\begin{enumerate}
\item ${w^*}^{\check{*}}=w$;
\item $\check{w}^{\check{*}*}=\check{w};$
\item $\check{N}(w^*)=N(w)^{*\B}$;
\item $w^{*\check{\sharp}}=w^{\sharp \check{*}}$;
\item $\check{w}^{\check{\sharp}*}=\check{w}^{\check{*}\sharp};$
\item $\langle w,\check{w}\rangle^{*_\B}=\langle \check{w}^{\check{*}}, w^*\rangle$;
\item $\langle w,w^*\rangle w^{* \check{\sharp}}= N(w)^{*_\B} w$
\end{enumerate}
for all $w\in \Pe$, $\check{w}\in \Pe^\vee$ [Ach1, 2.22].

We can now define
$$\widetilde{\J}=\A\oplus \Pe,$$
$$\tilde{1}=(1,0)\in H^0(X,\widetilde{\J}),$$
$$\begin{array}{r}
\widetilde{N}(a,w)=N_\B(a)+N(w)+\check{N}(w^*)-T_\B(a,\langle w,w^*\rangle)\\
=N_\B(a)+N(w)+N(w)^{*_\B}-T_\B(a,\langle w,w^*\rangle),
\end{array}$$
$$(a,w)^{\widetilde{\sharp}}=(a^{\sharp_\B}-\langle w,w^*\rangle,w^{*\check{\sharp}} -aw)$$
for $a\in \A$ and $w\in \Pe$. Then $(\widetilde{N},\widetilde{\sharp},\tilde{1})$ is a cubic form with adjoint
and base point on $\widetilde{\J}$ and
$$\widetilde{T}((a,w),(c,v))=T_\B(a,c)+T_\B(\langle w,v^*\rangle)+T_\B(\langle v,w^*\rangle)$$
for $a,c\in \A$ and $v,w\in \Pe$ is the trace. The induced Jordan algebra $\widetilde{\J}(\widetilde{N},\widetilde{\sharp},\tilde{1})$
is denoted by $\J(\B,\A,\Pe,N,*)$ [Ach1, 2.23] and called the {\it Tits process}.

In general, we do not get all possible cubic Jordan algebras over $X$ using the Tits process,
see [P-S-T2] for examples.

The Tits process $\J(\B,\A,\Pe,N,*)$ generalizes the classical Tits process over rings from [P-R1]: For $P\in X$,
$$\J(\B,\A,\Pe,N,*)_P\cong \J(\B_P,\A_P,\Pe_P,N_P,*_P)$$
where the right hand side is a classical Tits process over $\mathcal{O}_{P,X}$ in the sense of [P-R1].

\begin{example} Let $(N_\B,\sharp_\B, 1)$ be $\B$-admissible and let $(\A,\sts)$ be $\B$-ample.  $\B$ is canonically
a globally free left $\B$-module of rank one, denoted $_\B\B$. Hence $_\B\B$ is $\A$-admissible and the set of
$\A$-admissible pairs is identical with the set of pairs $(u*_\B,\beta N_\B)$ where $u\in H^0(X,\A^\times)$,
$\beta\in H^0(X,\sts')$ and $N_\B(u)=\beta\beta^{*_\B}$.
 We obtain the {\it classical} Tits process
$$\widetilde{\J}=\J(\B,\A,\,_\B\B,\beta N_\B,u*_\B)$$
with norm
$$\begin{array}{l}
\widetilde{N}(a,b)=N_\B(a)+\beta N_\B(b)+\beta^{-1}N_\B(ub^{*_\B})-T_\B(a,bub^{*_\B})=\\
N_\B(a)+\beta N_\B(b)+\beta^{*_\B}N_\B(b^{*_\B})-T_\B(a,bub^{*_\B}).
\end{array}$$
\end{example}

\begin{theorem} ([Ach1, 3.11]) Let $\J=\J(N_\J,\sharp_\J,1)$ be an Albert algebra over $X$ which contains a
subalgebra of the kind $\h(\B,*_\B)$, where $\B$ is an Azumaya algebra over $\mathcal{O}_X'$ of constant rank 9
together with an involution $*_\B$ and $\sts'$ an $\sts$-algebra of constant rank 2 with $\h(\B,*_\B)=\sts'$.
Then there exists an $\h(\B,*_\B)$-admissible left $\B$-module $\Pe$ of rank 1 and a $\h(\B,*_\B)$-admissible
pair $(N,*)$ such that the canonical embedding $\h(\B,*_\B)  \hookrightarrow \J$ can be extended to
an isomorphism $$\J(\B,\h(\B,*_\B),\Pe,N,*)\to \J.$$
\end{theorem}

\begin{example} Let $X$ be a proper integral algebraic scheme over a field $k$ of characteristic not 2 or 3
and let $k'$ be a separable quadratic field extension of $k$. Put $X'=X\times_k k'$ and $\sts'=\mathcal{O}_{X'}$.
\\ (i) With the above notation,
${\rm Pic}\,X'$ acts on $\check{H}^1(X',\B^\times)$ via
$(\lb,\Pe)\to \lb\otimes\Pe$ and $N_\B(\lb\otimes\Pe)\cong\lb^3\otimes N_\B(\Pe)$. In particular, if
$\lb\in {\rm Pic}\,X'$ has order 3, $ \lb\otimes\Pe$ is a left $\B$-module of rank one which admits a multiplicative
cubic norm (analogously as [P1, 2.10.]).
\\ (ii) Let $\Pe\in  \, _3 {\rm Pic} X'$ be a nontrivial line bundle with nondegenerate cubic form $N$. Let
$\star=\star_{\mathcal{O}_{X'}}$ be the  involution on $\mathcal{O}_{X'}$ induced by the nontrivial element in the Galois
group of $k'/k$.
The Tits process $\E=\J (\mathcal{O}_{X'},\sts, \mathcal{P}, N,\star)=\sts\oplus\Pe$ using $\Pe$ and $N$
is a commutative associative algebra over $X$ of rank 3 with norm
$$N((a,w))=a^3+N(w)+N(w)^{\star}$$ and trace
$$T((a,w),(c,v))=3ac+3\langle w,v^\star\rangle+3\langle v,w^\star\rangle$$
 for all $a,c\in \mathcal{O}_{X}$, $v,w\in\Pe$.
\end{example}

\subsection{The  \'etale Tits process} We rephrase the  \'etale Tits process ([P-T], [P-R2, p.~248])
 in the more general setting of locally ringed spaces: Let $\torus$ be a quadratic  \'etale algebra over $X$ with canonical involution $\can$.
Let $\E$ be a commutative associative $\sts$-algebra
 of constant rank 3 such that $\E^+=\J(N,\sharp,1)$.

 Specialize the Tits process to ${\sts}'=\torus$, the commutative
associative ${\sts}'$-algebra $\B=\E\otimes_{\sts}\torus$ and $*_\B=id_\E\otimes\can$. Then $\A=\h(\B,*_\B)
=\E\otimes 1$ can be identified with $\E$. We view $\B=\E_\torus$ as an algebra over $\torus$ to see that norm and adjoint of $\E$ canonically
extend to $\B$.
(If $\torus=k'\otimes \sts$ is a quadratic  \'etale algebra which is
defined over $k$ with $k'$ a cubic field extension, the residue class algebras are cubic  \'etale over $k'(P)$.)
 $\J=\J(\B,\h(\B,*_\B),\Pe,N,*)$ is called an {\it  \'etale Tits process}.
Since for all $P\in X$, the residue class algebra $\J(P)$ is an  \'etale Tits process over the residue class
field $k(P)$, $\J(P)$ is an absolutely simple Jordan algebra of degree 3 and dimension 9 for all $P\in X$. Thus $\J(P)$
is the symmetric elements of a central simple associative algebra of degree 3 with involution of the second kind
for all $P\in X$ [P-T, p.~91].
This implies that $\J$ does not decompose into the direct sum of two non-trivial ideals
$\J_1\oplus\J_2$.

It remains to check if $\J$ itself perhaps
is the symmetric elements of an Azumaya  algebra over ${\sts}'$ of rank 9 with involution of the second kind.

Every  homomorphism $f$ between two  \'etale Tits processes preserving norms and units
is an isomorphism.
($f(P)$ preserves norms and units, hence $f(P)$ must be an isomorphism by [P-T, p.~91].
 Thus so must be $f$ by Nakayama's Lemma.)

\subsection{Algebras admitting an involution}

\begin{proposition} Let $K$ be a separable quadratic field extension of $k$ with ${\rm Gal}(K/k)=\langle\omega\rangle$.
Let $X'=X\times_k K$ and let $\E$ be a vector bundle over $X'$ of rank $r$. The involution $\omega$ on $K$ induces
an involution on $X'$ we also call $\omega$. If $h:\E\to^\omega\E^\vee$ is a nondegenerate hermitian form (i.e., $h$ is an
 isomorphism
such that $h=\,^\omega h^\vee$), then $\B=\E nd_{X'} (\E)$ is an Azumaya algebra of rank $r^2$ over $X'$ which admits an
involution $*_\B$ such that $*_\B|_{\mathcal{O}_{X'}}=\omega$ and $\h (\mathcal{O}_{X'},*_\B)=\sts$.
\end{proposition}

\begin{proof} As in [KMRT, p.~43] define $$\sigma_h(f)=h^{-1}\circ\,^\omega f^t \circ h$$
 where $^\omega f^t:\,^\omega \E^t\to \,^\omega \E^t$
is the transpose of $f\in\E nd_{X'}(\E)$. Then $\sigma_h$ is an involution on $\B$ (the {\it adjoint involution})
such that $\sigma_h(a)=\omega(a)$ for all $a\in\mathcal{O}_{X'}$ and
$$h(x,f(y))=h(\sigma_h(f)(x),y)$$ for all $x,y\in\E$.
\end{proof}

Over  a field $k$, the map $h\to\sigma_h$ defines a $1$-$1$-correspondence between nondegenerate hermitian forms
$h: E\to^\omega E^\vee$ and
 involutions $\sigma$ on ${\rm End}_K(E)$ of the second kind such that $\sigma(a)=\omega(a)$ for all $a\in  k$
 [KMRT, I.(4.2)].

\begin{remark} Let $X$ be a proper scheme over a perfect field of characteristic not two and let
$K$ be a quadratic field extension of $k$ with ${\rm Gal}(K/k)=\langle\omega\rangle$.
Let $X'=X\times_k K$ and let $\E$ be a vector bundle over $X'$ of rank $3$. Suppose that
 $h:\E\to^\omega\E^\vee$ is a nondegenerate hermitian form on $\E$. We use the results from [AEJ1]:
\\(1)  If $\E$ is absolutely indecomposable then $h$ is uniquely determined up to isomorphism and  an
 invertible scalar factor.
 \\(2) If $\E\cong\lb\oplus\F$ for a line bundle $\lb$ and an absolutely indecomposable
 vector bundle $\F$ of rank 2 then $^\omega\lb\cong\lb^\vee$, $^\omega\F\cong\F^\vee$ by the Theorem for Krull-Schmidt.
 Both
 $\lb$ and $\F$ must admit a nondegenerate hermitian form $h_0$ respectively $h_1$ which again
 is uniquely determined up to isomorphism and an invertible scalar factor. Thus $h\cong \mu_0 h_0\perp \mu_1 h_1$ for suitable
 invertible $\mu_i\in H^0(X,\sts)$.
\\(3) If $\E\cong\lb\oplus\M\oplus \N$ for  line bundles $\lb$, $\M$ and $\N$,
  then either $^\omega\lb\cong\lb^\vee$, $^\omega\M\cong\M^\vee$ and $^\omega\N\cong\N^\vee$
  and all three bundles carry a nondegenerate hermitian form which each is uniquely determined up to isomorphism and  an
 invertible scalar factor, or $(\E,h)\cong (\lb,h_0)\perp\\
  (\M\oplus ^\omega\M^\vee,\mathbb{H})$ where $\mathbb{H}$
 denotes a hyperbolic space of rank 2.
\\(4) Let $\E$ be indecomposable and $\E=tr_{L/K}(\N)$ for some line bundle $\N$ over $Y=X'\times_KL$,
 where $L/K$ is a field
extension of degree 3. Suppose that there is given an involution on $L$ extending $\omega$ (which we also
call $\omega$). Then every nondegenerate hermitian form $h_0:\N\to\,^\omega\N^\vee$ induces a nondegenerate hermitian form $h:
tr_{L/K}(\N)\to ^\omega tr_{L/K}(\N)^\vee$ via $h=tr_{L/K}(h_0)$ [AEJ1, p.~1349].
 \\(5) Let $\E\cong\lb\oplus\F$ for a line bundle $\lb$ and an  indecomposable
 vector bundle $\F=tr_{L/K}(\N)$, with some line bundle $\N$ over $Y=X'\times_KL$, where $L/K$ is a field
extension of degree 2. Then $^\omega\lb\cong\lb^\vee$, $^\omega\F\cong\F^\vee$ by the Theorem for Krull-Schmidt
 and both  $\lb$ and $\F$ must admit a nondegenerate hermitian form $h_0$ respectively $h_1$.
Assuming that there is an involution $\omega$ on $L$ extending $\omega$, every nondegenerate hermitian form $h_0:\N\to\,^\omega\N^\vee$ again induces a nondegenerate hermitian form $h:
tr_{L/K}(\N)\to ^\omega tr_{L/K}(\N)^\vee$ via $h=tr_{L/K}(h_0)$.
\end{remark}

\section{The first Tits construction}

 Let $(X,\sts)$ be a locally ringed space and let $\A$ be a unital associative $\sts$-algebra, $(N_\A,\sharp_\A,1)$ a cubic form with
adjoint and base point on $\A$, $\Pe\in {\rm Pic}_l\,\A$ such that $N_\A(\Pe)\cong\sts$ and $N$ a cubic
norm on $\Pe$. Let $\A^+=J(N_\A,\sharp_\A,1)$ and $N_\A(xy)=N_\A(x)N_\A(y)$ for all $x,y\in\A$. Define
$$\begin{array}{l}
\widetilde{\J}=\A\oplus \Pe\oplus\Pe^\vee,\\
\widetilde{1}=(1,0,0)\in H^0(X,\J),\\
\widetilde{N}(a,w,\check{w})=N_\A(a)+N(w)+\check{N}(\check{w})-T_\A(a,\langle w,\check{w}\rangle)\\
(a,w,\check{w})^{\widetilde{\sharp}}=(a^{\sharp_\A}-\langle w,\check{w}\rangle, \check{w}^{\check{\sharp}}-
aw, w^\sharp-\check{w}a)
\end{array}$$
for $a\in\A$, $w\in\Pe$, $\check{w}\in\Pe^\vee$. Then $(\widetilde{N},\widetilde{\sharp},\widetilde{1})$
is a cubic form with adjoint and base point on $\widetilde{\J}$ and with trace form
$$\widetilde{T}((a,w,\check{w}),(c,v,\check{v}))=T_\A(a,c)+T_\A(\langle w, \check{v}\rangle)+T_\A(\langle v,
\check{w}\rangle).$$
The induced Jordan algebra $\J(\widetilde{N},\widetilde{\sharp},\widetilde{1})$ is denoted by $\J(\A,\Pe,N)$
and called a {\it first Tits construction} [Ach1, 2]. $\A^+$ identifies canonically with a subalgebra of
$\J(\A,\Pe,N)$.\\
If
$$\begin{array}{l} \sts'=\sts\oplus\sts, \, (r,s)^{*_\B}=(s,r),\\
\B=\A\oplus\A^{\rm op}, \, (a,c)^{*_\B}=(c,a),\, 1_\B=(1,1),\\
N_\B(a,c)=(N_\A(a),N_\A(c)), \, (a,c)^{\sharp_\B}=(a^{\sharp_\A},c^{\sharp_\A}),\\
\sts^0=\{(r,r)\,|\,r\in\sts\},\, \A^0=\{(a,a)\,|\,a\in\A\},\\
\Pe^0=\Pe\oplus\Pe^\vee,\\
N^0(w,\check{w})=(N(w),\check{N}(\check{w})),\, (w,\check{w}^{*^0})=(\check{w},w),
\end{array}$$
then
$$\J(\B,\A^0,\Pe^0,N^0,*^0)\cong \J(\A,\Pe,N)$$ [Ach1, 2.25].

For $P\in X$, $$\widetilde{\J}(\A,\Pe,N)_P\cong \widetilde{\J}(\A_P,\Pe_P,N_P)$$
where the right hand side is a  classical first Tits construction over $\mathcal{O}_{P,X}$ in the sense of [P-R1].

\begin{example}
This first Tits construction generalizes the classical first Tits construction from [P-R1]:
If we choose $\Pe$ to be the left $\A$-module $_\A\A$ itself then the norms on $_\A\A$ are exactly the cubic
maps $\beta N_\A$ for $\beta\in H^0(X,\sts^\times)$ and
$$\widetilde{\J}=\J(\A,_\A\A,\beta N_\A)$$
is a classical first Tits construction with norm
$$\widetilde{N}(a_0,a_1,a_2)=N_\A(a_0)+\beta N_\A(a_1)+\beta^{-1}N_\A(a_2)-T_\A(a_0,a_1a_2).$$
\end{example}

\begin{theorem} ([Ach1, 3.9]) If $\J=\J(N_\J,\sharp_\J,1)$ is an Albert algebra over $X$ containing a subalgebra of the kind $\A^+$ with $\A$ an
Azumaya algebra over $X$, then there exist suitable $\Pe$ and $N$ such that $\J(\A,\Pe,N)\cong \J$.
\end{theorem}

\begin{remark} Albert algebras over integral schemes  were investigated by
Parimala-Suresh-Thakur [P-S-T1, 2].
In [P-S-T2],  generalized first and second Tits constructions for Albert algebras
over a domain $R$ such that $2,3\in R^\times$ were introduced.
 Over an integral scheme over a field $k$ of characteristic not 2 or 3, the first Tits construction coincides with the first Tits construction
  given above  when starting
with an Azumaya algebra of constant rank 9, the second with the Tits process of Achhammer starting with
the subsheaf of hermitian elements of an Azumaya algebra $\mathcal{B}$ of constant rank 9 over $\mathcal{O}_X'$ with involution
$*_\mathcal{B}$ such that $*_\mathcal{B}|_{\mathcal{O}_X'}=*$, where $(X,\mathcal{O}_X')$ is a locally ringed
space with involution $*$ such that  $\mathcal{H}(\mathcal{O}_X',*)=\mathcal{O}_X$.
\end{remark}

\begin{theorem} ([Ach1, 2.28]) Suppose $\B$, $\mathcal{O}_X'$, $*_\B$, $(N_\B,\sharp_\B,1)$, $\Pe$ and $(N,*)$ are as in
3.3 so that we have the Tits process $\J(\B,\A,\Pe,N,*)$. In particular,
$\B$, $\mathcal{O}_X'$, $(N_\B,\sharp_\B,1)$, $\Pe$ and $N$ satisfy the assumptions for the general
 first Tits construction $\J(\B,\Pe,N)$. Define an involution
 $$*:\J(\B,\Pe,N)\to \J(\B,\Pe,N)$$ via
 $$(b,w,\check{w})^*=(b^{*_\B},\check{w}^{\check{*}},w^*)$$
 for $b\in \B$, $w\in \Pe$, $\check{w}\in \check{P}$. Then
 $$\J(\B,\A,\Pe,N,*)\cong \mathcal{K}\subset \h (\J(\B,\Pe,N),*),$$
 i.e., the algebra $\J(\B,\A,\Pe,N,*)$ is isomorphic to a subsheaf $\mathcal{K}$ of the $\h (\mathcal{O}_X',*_\B)$-Jordan
 algebra $\h(\J(\B,\Pe,N),*)$. For $\A=\h (\B,*_\B)$ we have
 $$\mathcal{K}= \h (\J(\B,\Pe,N),*).$$
\end{theorem}

Hence the Tits process can be embedded into a  general first Tits construction, analogously as shown
in [P-R1], [M1].

\begin{lemma} ([Ach1, 3.10]) Let $\J=\J(N_\J,\sharp_\J,1)$ be an Albert algebra over $X$ containing a subalgebra of the kind
$\h(\B,*_\B)$ where $\B$ is an Azumaya algebra over $\sts'$ of constant rank 9 with involution $*_\B$ and $\sts'$ an $\sts$-algebra
of constant rank 2 with $\h(\sts',*_\B)=\sts$. Then $\B^+\cong\h(\B,*_\B)\otimes\sts'$ is an $\sts'$-sub Jordan
algebra of $\J'=\J\otimes\sts'$.
\end{lemma}

\begin{proposition} Assume the situation of Lemma 4. Then $\J'\cong\J(\B,\Pe_0,N_0)$
for suitable $\Pe_0$ and $N_0$.
\end{proposition}

The proof follows immediately from Lemma 4 and Theorem 3.
This observation restricts the $\sts$-module structure for these Albert algebras over proper schemes:

\begin{corollary} Let $X$ be a proper scheme over a perfect base field $k$, $k'/k$ a quadratic
field extension with ${\rm Gal}(k'/k)=\langle\omega\rangle$, $X'=X\times_kk'$ and $\sts'=\mathcal{O}_{X'}$.
Let $\J$ be an Albert algebra as in Lemma 4, $\Pe_0$ and $N$ as in Proposition 3.
\\(i) If $\Pe_0$ is indecomposable over $X'$ and there is an indecomposable $\sts$-module $\F$ of rank 9 such that
 $\F\otimes \mathcal{O}_{X'}\cong\Pe_0$, then
 $$\Pe\cong \F\oplus \F^\vee$$
 as $\sts$-module.
\\ (ii) If $\Pe_0$  is indecomposable over $X'$ and not defined over $X$, then $ ^\omega\Pe_0\cong\Pe_0^\vee$
and
$$\Pe\cong tr_{k'/k}(\Pe_0)$$
is an indecomposable $\sts$-module of rank 18.
\\ (iii) If the Krull-Schmidt decomposition of $\Pe_0$ is given by
 $\Pe_0\cong \F_1\oplus\Pe_1$ with $\F_1$ defined over $X$ and $\Pe_1$ not defined over $X$, then
 $\Pe_1^\vee\cong^\omega\Pe_1$ and
 $$\Pe\cong \F_1\oplus \F_1^\vee\oplus tr_{k'/k}(\Pe_1)$$
with $tr_{k'/k}(\Pe_1)$ an indecomposable $\sts$-module of rank $2\cdot {\rm rank}\,\Pe_1$.
\\ (iv) If the Krull-Schmidt decomposition of $\Pe_0$ is given by
 $\Pe_0\cong \F_1\oplus\Pe_1$ with $\F_1$ and $\Pe_1$ not defined over $X$ and of different ranks, then
 $\F_1^\vee\cong^\omega\F_1$, $\Pe_1^\vee\cong^\omega\Pe_1$ and
 $$\Pe\cong tr_{k'/k}(\F_1)\oplus tr_{k'/k}(\Pe_1)$$
is the direct sum of two indecomposable $\sts$-modules.
\end{corollary}

\begin{proof} Since $\h(\B,*_\B)$  is contained in $\J$, there are also suitable $\Pe$ and $(N,*)$ such that
$\J\cong \J(\B,\h(\B,*_\B),\Pe,N,*)$ and $\J'\cong\J(\B,\Pe_0,N_0)$
for suitable $\Pe_0$ and $N_0$ by Proposition 3. We can use the Theorem of
Krull-Schmidt for vector bundles and the theory in [AEJ1] to conclude that
$$\Pe\otimes_{\sts }\mathcal{O}_{X'}\cong \Pe_0\oplus \Pe_0^\vee$$
as $\mathcal{O}_{X'}$-modules, where ${\rm rank}_ {\mathcal{O}_{X'}}\Pe_0=9$.\\
If $\Pe_0$ is indecomposable over $X'$, then either
there is an indecomposable $\sts$-module $\F$ of rank 9 such that
$$\Pe\cong \F\oplus \F^\vee \text{ and }\F\otimes \mathcal{O}_{X'}\cong\Pe_0$$
or $\Pe_0$ is not defined over $X$ and
$$\Pe\cong tr_{k'/k}(\Pe_0)  \text{ with } ^\omega\Pe_0\cong\Pe_0^\vee$$
is an indecomposable vector bundle of rank 18.
This settles (i) and (ii).\\
Let
$$\Pe_0=\F_1\oplus \dots\oplus\F_s\oplus\dots\oplus \Pe_m,$$
$0\leq s, m\leq 9$, be the Krull-Schmidt decomposition of $\Pe_0$ into
indecomposable sumands. Suppose w.l.o.g. that the first $s$, $0\leq m\leq s$ of these are defined over $X$.
A case-by-case study implies the assumptions (iii) and (iv).
\end{proof}

Note that the method used in the proof of the above corollary can be continued to consider more elaborate
cases of Krull-Schmidt decompositions of $\Pe_0$ if desired.

\begin{corollary} Let $X$ be a scheme over the affine scheme $Y={\rm Spec}\,k$ and suppose $H^0(X,\mathcal{O}_X)=k$.
 Let $k'/k$ be a separable quadratic field extension, put $X'=X\times_kk'$.
Let $\J=\J(N_\J,\sharp_\J,1)$ be an Albert algebra over $X$ containing a subalgebra of the kind
$\h(\B,*_\B)$ where $\B$ is an Azumaya algebra over $\mathcal{O}_{X'}$ of constant rank 9 with involution $*_\B$
such that $\h(\mathcal{O}_{X'},*_\B)=\sts$. If $\J'$ is defined over $k'$ then $\J$ is defined over $k$.
\end{corollary}

\begin{proof} If $\J'$ is defined over $k'$ then it is globally free as $\mathcal{O}_{X'}$-module and thus so is
$\J$, implying that it must be defined over $k$.
\end{proof}

\begin{example} In the situation of 3.4 let $\E=\sts\times \sts\times\sts$. Then,
by [P-T, 4.9],  the  \'etale Tits process satisfies $\J(P)\cong \h_3(\torus(P),g)$ for all $P\in X$,
 with $g$ depending on $P$. It would be interesting to know how $\J$ looks like globally.

Let  $k'$ be a quadratic  \'etale field extension of $k$, put $X'=X\times_k k'$. Then $\B=\E\otimes \torus=
\E\otimes \mathcal{O}_{X'}\cong  \mathcal{O}_{X'} \times \mathcal{O}_{X'}  \times \mathcal{O}_{X'}$.
 Since $\check{H}^1(X',\B^\times)={\rm Pic}\,X'\oplus {\rm Pic}\,X'\oplus {\rm Pic}\,X'$,
 every left $\B$-module $\Pe$
of rank one satisfies $\Pe\cong \lb\oplus\M\oplus \s$ with invertible $ \mathcal{O}_{X'}$-modules $\lb$, $\M$ and $\s$.
Thus $ N_\B(\Pe)\cong  \mathcal{O}_{X'} $ iff $\lb\otimes \M\otimes \s\cong  \mathcal{O}_{X'}$.

Suppose that $X$ is an integral scheme over $k$ such that $2\in H^0(X,\sts)$ is invertible. $\J$  contains
$\mathcal{E}^+=\sts^+\times\sts^+\times\sts^+$. $H^0(X,\mathcal{E})$ has a complete orthogonal system
of primitive idempotents $(c_1,c_2,c_3)$.
Hence, by Remark 3, there are line bundles $\J_{i,j}$ such that
 $\J_{i,i}=\sts c_i$ ($1\leq i\leq 3$),
 $$\J_{i,j}=\J_{j,i}=\{x\in \J\,|\,xc_i=\frac{1}{2}x=xc_j \}$$
($1\leq i<j\leq 3$)
and
$$\J=\J_{1,1}\oplus \J_{1,2}\oplus \J_{1,3}\oplus \J_{2,1}\oplus \J_{2,2}\oplus \J_{2,3}\oplus \oplus \J_{3,1}\oplus
\J_{3,2}\oplus \J_{3,3}$$
as an $\sts$-module.
\end{example}

\begin{example} Let $X$ be a scheme having $2,3\in R:=H^0(X,\mathcal{O}_X) $.
 Every first Tits construction  over $X$ starting with $\mathcal{O}_X$ has the form
 $\J (\mathcal{O}_X, \mathcal{L}, N)$ for some line bundle
 $\mathcal{L} \in \, _3 {\rm Pic}\,X$ and some cubic norm $N$ on $\mathcal{L}$. We also note that
  $\bigwedge^3(\J (\mathcal{O}_X, \mathcal{L},  N))\cong\sts$.
 More precisely:

Let $\lb^\vee= \mathcal{H}om_{\sts}(\lb,\sts)$ and $\langle w,\check{w}\rangle=\check{w}(w)$ be the canonical pairing $\lb
\times\lb^\vee\to \sts.$ There exists a uniquely determined cubic norm $\check{N}:\lb^\vee\to \sts$ and
uniquely determined adjoints $\sharp:\lb\to \lb^\vee$ and $\check{\sharp}:\lb^\vee\to \lb$ such that
 $\langle w,w^\sharp\rangle=N(w)1$; $\langle \check{w}^{\check{\sharp}},\check{w}\rangle=\check{N}(\check{w})1$
and $w^{\sharp \, \check{\sharp}}=N(w)w$
for $w$ in $\lb$, $\check{w}$ in $\lb^\vee$.
 The norm, adjoint and trace of
 $$\J (\mathcal{O}_X, \mathcal{L}, N)=\sts\oplus \lb\oplus\lb^\vee$$
  are given by
$$\begin{array}{l}
\widetilde{N}(a,w,\check{w})=a^3+N(w)+\check{N}(\check{w})-3a\langle w,\check{w}\rangle\\
(a,w,\check{w})^{\widetilde{\sharp}}=(a^2-\langle w,\check{w}\rangle, \check{w}^{\check{\sharp}}-aw, w^\sharp-\check{w}a)\\
\widetilde{T}((a,w,\check{w}),(c,v,\check{v}))=3ac+3\langle w, \check{v}\rangle+3\langle v,
\check{w}\rangle
\end{array}$$
for $a,c\in\sts$, $v,w\in\lb$, $\check{v},\check{w}\in\lb^\vee$.
This construction yields examples of algebras $\A=\J (\sts, \lb,  N)$, where
 $\A(P)$ is a classical first Tits construction starting with $k(P)$ and thus a
  cubic  \'{e}tale algebra over $k(P)$ for all $P\in {\rm Spec}\,X$.
If $X$ is an integral scheme, $\A$ is commutative associative, since so is $\A(P)$ for all $P\in X$.
\end{example}

\subsection{The \'etale first Tits construction}
Let $\A$ be a commutative associative $\sts$-algebra of constant rank 3 such that $\A^+=\mathcal{J}(N_\A,\sharp_\A,1)$
and where $\A(P)$ is a cubic
\'etale $k(P)$-algebra for all $P\in X$. The first Tits construction $\J=\mathcal{J}(\A,\Pe,N)$ is called an {\it
\'etale first Tits construction}.

 Use the notation of 3.4 where the the \'etale Tits process was formulated:
Let $\torus$ be a quadratic \'etale algebra over $X$ with canonical involution $\can$ and let $\E$ be a commutative associative $\sts$-algebra
 of constant rank 3 such that $\E^+=\mathcal{J}(N,\sharp,1)$.
 Specialize the Tits process to ${\sts}'=\torus$, the commutative
associative ${\sts}'$-algebra $\B=\E\otimes_{\sts}\torus$ and $*_\B=id_\E\otimes\can$.
If $\torus=\sts\times\sts$ is a split quadratic \'etale algebra with exchange involution $\can$ then
$\B=\E\oplus\E$ with exchange involution and $H^1(X,\B)=H^1(X,\E)\oplus H^1(X,\E)$.
Hence any $\Pe\in H^1(X,\B)$ is of the type $\Pe\cong\Pe_1\oplus\Pe_2$ with $\Pe_i\in H^1(X,\E)$.
 A tedious calculation shows that, as over
fields, in this case the \'etale Tits process
$\mathcal{J}(\B,\A,\Pe,N,*)\cong \E\oplus \Pe_1\oplus \Pe_1^\vee$
with $\Pe_1\in H^1(X,\E)={\rm Pic}_l\,\E$ becomes an \'etale first Tits construction starting with $\E$.

\begin{example}  ([Ach1, 4.2]) Let $\A=\sts\times\sts\times\sts$. Then $\A^+=\J(N_\A,\sharp_\A,1)$ by Example 1 (i).
Since $\check{H}^1(X,\A^\times)={\rm Pic}\,X\oplus {\rm Pic}\,X\oplus {\rm Pic}\,X$, every left $\A$-module $\Pe$
of rank one satisfies $\Pe\cong \lb\oplus\M\oplus \s$ with invertible $\sts$-modules $\lb$, $\M$ and $\s$.
Thus $ N_\A(\Pe)\cong \sts$ iff $\lb\otimes \M\otimes \s\cong \sts$.
Choose an isomorphism $\alpha:\lb\otimes \M\otimes \s\to\sts$, then $N(x,y,z)=\alpha(x\otimes y\otimes z)$
defines a norm on $\Pe$ and the  \'etale first Tits construction
 $\J=\J(\A,\Pe,N)$ is a Jordan algebra over $X$ with
$$\begin{array}{l}
N_J(a,w,\check{w})=N_\A(a)+N(w)+\check{N}(\check{w})-T_\A(a,\langle w,\check{w}\rangle),\\
(a,w,\check{w})^{\sharp_\J}=(a^{\sharp_\A}-\langle w,\check{w}\rangle,\check{w}^{\widetilde{\sharp}}-aw, w^\sharp-\check{w}a)
\end{array}$$
and
$$\begin{array}{l}
(x,y,z)^\sharp=(y\otimes z,z\otimes x, x\otimes y),\\
(\check{x},\check{y},\check{z})^{\widetilde{\sharp}}=(\check{y}\otimes \check{z},\check{z}\otimes \check{x}, \check{x}
\otimes\check{y}),\\
\check{N}(\check{x},\check{y},\check{z})= \check{x}\otimes \check{y}\otimes \check{z},\\
\langle (x,y,z),(\check{x},\check{y},\check{z})\rangle=(\langle x,\check{x}\rangle,\langle y,
\check{y}\rangle,\langle z,\check{z}\rangle).
\end{array}$$
for $(x,y,z)$ in $\Pe$, $(\check{x},\check{y},\check{z})$ in $\Pe^\vee$.
(We canonically identify $\lb\otimes\M\cong \check{\s}$ etc.)

 $J(\A,\Pe,N)$ has the direct sum of line bundles
 $ (\sts\oplus\sts\oplus\sts)\oplus ( \lb\oplus\M\oplus \s)\oplus
(\check{ \lb}\oplus\check{\M}\oplus \check{\s})$ as underlying module structure.
\end{example}

For certain \'etale algebras over proper schemes we can restrict their $\sts$-module structure:

\begin{lemma} Let $X$ be a proper scheme over a perfect field $k$.
\\ (i) Let $k'/k$ be a  quadratic field extension of $k$ with ${\rm Gal}(k'/k)=\langle \omega\rangle$.
Let $\torus=k'\otimes\sts $.
For $\A=\sts\times \torus$, every left $\A$-module $\Pe$
of rank one satisfies $\Pe\cong \lb\oplus\F$ with $\F\in {\rm Pic}\torus$ and an invertible $\sts$-module $\lb$.
Every left $\A$-module $\Pe$ of rank one with $N_\A(\Pe)\cong\sts$ is one of the following:\\
(a)  a direct sum of line bundles $\lb_0,\M_0,\s_0$ over $X$ with
$\lb_0\otimes\M_0\otimes\s_0\cong\mathcal{O}_{X}$; i.e.,
 $$\Pe\cong \lb_0\oplus\M_0\oplus \s_0;$$
(b) there is a line bundle $\s_0$ over $X$ and a line bundle $\lb$  over $X'$, which is not defined over $X$,
satisfying
$\lb\otimes\, ^\omega\lb\otimes\s_0'\cong\mathcal{O}_{X'}$
such that $tr_{k'/k}(\lb)\in {\rm Pic} (\torus)$ and
$$\Pe\cong tr_{k'/k}(\lb)\oplus\s_0.$$
\noindent (ii)  Let $k'/k$ be a separable cubic field extension of $k$. Let $k^+$ be a Galois
extension of $k$ containing $k'$. For $\A=k'\otimes \sts$, every left
 $\A$-module $\Pe$ of rank one with $N_\A(\Pe)\cong\sts$ is one of the following:\\
 (a) a direct sum of line bundles $\lb_0,\M_0,\s_0$ over $X$ with
$\lb_0\otimes\M_0\otimes\s_0\cong\mathcal{O}_{X}$; i.e.,
 $$\Pe\cong \lb_0\oplus\M_0\oplus \s_0;$$
(b) there is a line bundle $\lb$  over $X'$ which is not defined over $X$ whose $k'/k$-conjugates
satisfy
$\lb^+\otimes\, ^{\sigma_1}\lb^+\otimes \, ^{\sigma_2}\lb^+\cong\mathcal{O}_{X^+}$
such that
$$\Pe\cong tr_{k'/k}(\lb).$$
\end{lemma}

\begin{proof} (i) Obviously, each such $\A$ is defined over $k$,
thus there is a cubic form with adjoint and base point on $\A$ such that $\A^+=\J(N_\A,\sharp_\A,1)$.
Put $X'=X\times_k k'$ and let $\J=J(\A,\Pe,N)$. Then
$$\J\otimes_{\sts} \mathcal{O}_{X'}=J(\A\otimes \mathcal{O}_{X'},\Pe\otimes\mathcal{O}_{X'},N\otimes\mathcal{O}_{X'})$$
and $\A\otimes \mathcal{O}_{X'}\cong \mathcal{O}_{X'}\times\mathcal{O}_{X'}\times\mathcal{O}_{X'}$.
Therefore there exist $\lb,\M,\s\in {\rm Pic}X'$ with $\lb\otimes\M\otimes\s\cong\mathcal{O}_{X'}$, such that
$$\Pe\otimes\mathcal{O}_{X'}\cong \lb\oplus\M\oplus \s.$$
This restricts the $\sts$-module structure of $\Pe$ and yields the assertion using [AEJ1].
By the previous example, every left $\A$-module $\Pe$
of rank one satisfies $\Pe\cong \lb\oplus\F$ with $\F\in {\rm Pic}\torus$ and an invertible $\sts$-module $\lb$.
 We have $ N_\A(\Pe)\cong \sts$ iff $\lb\otimes N_\torus(\F)\cong \sts$ iff $\lb^\vee\cong N_\torus(\F)$.
\\ (ii) is proved analogously.
\end{proof}

\subsection{The first Tits construction starting with an algebra of rank 5}

\begin{lemma}  Let $\D$ be a quaternion algebra over $X$ and $\A=\sts\times\D$.
 \\ (i) Every left $\A$-module $\Pe$ of rank one and norm one satisfies
$$\Pe\cong  \lb\oplus \Pe_0$$
where $\lb$ is a line bundle over $X$ and $\Pe_0\in {\rm Pic}_l\,\D$ has $N_\D(\Pe_0)\cong \lb^\vee.$
 Hence the underlying module structure for a first Tits construction $\J=J(\A,\Pe,N)$ is given by
$$\sts\oplus \D\oplus \lb\oplus \Pe_0\oplus \lb^\vee\oplus \Pe_0^\vee.$$
\\ (ii) Let $\D=\E nd_X(\E)$. Then every left $\A$-module $\Pe$ of rank one and norm one satisfies
$$\Pe\cong  ({\rm det}\,\E)^\vee \otimes {\rm det}\,\F\oplus  \E\otimes\F^\vee$$
where $\F$ is a locally free $\sts$-module of constant rank 2.
 Hence the underlying module structure for a first Tits construction $\J=J(\A,\Pe,N)$ is given by
$$\sts\oplus \E nd_X(\E)\oplus [({\rm det}\,\E)^\vee \otimes {\rm det}\,\F\oplus  \E\otimes\F^\vee]
\oplus [({\rm det}\,\E) \otimes ({\rm det}\,\F)^\vee \oplus  \F\otimes\E^\vee]. $$
\end{lemma}

\begin{proof} $\A^+=J(N_\A,\sharp_\A,1)$ by Example 1.
\\ (i) We have $\check{H}^1(X,\A^\times)={\rm Pic}\, X\oplus {\rm Pic}_l\,\D $, thus  every left $\A$-module
$\Pe$ of rank one satisfies
$\Pe=\lb\oplus \Pe_0$ for some $\Pe_0\in {\rm Pic}_l\,\D$ and some invertible $\sts$-module $\lb$.
 Now $N_\A(\Pe)\cong \sts$ if and only if
$\lb\otimes N_\D(\Pe_0)\cong \sts$
if and only if $N_\D(\Pe_0)\cong \lb^\vee.$
\\ (ii) Each $\Pe_0\in {\rm Pic}_l\,\D$ is of the form $\Pe_0\cong \E\otimes\F^\vee$,
where $\F$ is a locally free $\sts$-module of constant
rank 2 , and $ N_\D(\Pe_0)\cong {\rm det}\,\E\otimes ({\rm det}\,\F)^\vee$ [P1].
Therefore every locally free left $\A$-module of rank 1 which admits a multiplicative cubic form is isomorphic to
$$\Pe\cong  ({\rm det}\,\E)^\vee \otimes {\rm det}\,\F\oplus  \E\otimes\F^\vee.$$
\end{proof}

\subsection{Albert algebras arising from the first Tits construction}
Let $\A$ be an Azumaya algebra over $X$. Then $\J(\A,\Pe,N)$ is an Albert algebra over $X$.

\begin{example} ([Ach1, 4.2]) For every locally free $\sts$-module $\mathcal{E}$
 of constant
rank 3, $\A=\E nd_X(\E)$ is an Azumaya algebra of rank 9 and $\A^+=\J(N_\A,\sharp_\A,1)$ with $N_\A$ the usual determinant,
$\sharp_\A$ the usual adjoint. The locally free left $\A$-modules all have the form
 $$\Pe=\E\otimes\F^\vee=\h om_X(\F,\E),$$
 $\F$ being another locally free $\sts$-module of constant rank 3. Since $N_\A(\Pe)\cong {\rm det}\, \E\otimes
{\rm det}\, \F^\vee$, we have $N_\A(\Pe)\cong\sts$ iff there exists an isomorphism $\alpha: {\rm det}\, \E\to
{\rm det}\, \F$. Fixing such an isomorphism, there exists a unique form $N:\Pe\to\sts$ with
$$(\alpha\circ{\rm det})(g)=N(g)id_{{\rm det}\F}$$
for $g$ in $\Pe$. This $N$ is a cubic norm on $\Pe$. Moreover,
$$\begin{array}{l}
\Pe^\vee=\E^\vee\otimes \F=\h om_X(\E,\F)\\
\langle g,f\rangle=g\circ f \text{ for } g\text{ in } \Pe, f \text{ in } \Pe^\vee.
\end{array}$$
The adjoint $\sharp:\Pe\to\Pe^\vee$ of $N$ satisfies
$$g\circ g^\sharp=N(g)id_\E.$$
The first Tits construction $\J(\A,\Pe,N)$ is an Albert algebra over $X$.
\end{example}

\begin{proposition} If $\comp$ is a split octonion algebra, then there are  suitable $\Pe$ and $N$ such that
$\mathcal{H}_3(\comp)$ is a first Tits construction and so $\mathcal{H}_3(\comp)=\J({\rm Mat}_3(\mathcal{O}_X), \Pe,N)$.
 If the Krull-Schmidt Theorem
holds over $X$ and $\comp={\rm Zor}\,(\torus,\alpha)$, then $\Pe\cong  {\F^\vee}^3$ as $\sts$-module
for some locally free
$\sts$-module $\F$ of constant rank 3 with trivial determinant and
$\F^3\oplus{\F^\vee}^3\cong\torus^3\oplus { \torus^\vee}^3$ as $\sts$-module.
\end{proposition}

\begin{proof} If $\comp_1$ is a composition subalgebra of $\comp_2$, then $\mathcal{H}_3(\comp_1,\Gamma)\subset
\mathcal{H}_3(\comp_2,\Gamma)$
is a subalgebra. Thus, if $\comp$ is a split composition algebra, then $\mathcal{H}_3(\comp)$ contains an isomorphic
copy of ${\rm Mat}_3(\mathcal{O}_X)^+$.  The assertion now follows from Example 9, since
${\rm Mat}_3(\mathcal{O}_X)=\E nd_X(\sts\oplus\sts\oplus\sts)$.
A left ${\rm Mat}_3(\mathcal{O}_X)$-module $\Pe$ has the form  $\Pe\cong \sts^3\otimes \F^\vee$ for some locally free
$\sts$-module $\F$ of constant rank 3 with trivial determinant.
Thus $\Pe\cong  {\F^\vee}^3$ as $\sts$-module. If the Krull-Schmidt Theorem holds over $X$
and  $\comp={\rm Zor}\,(\torus,\alpha)$ (for the definition, see [P1]), then
$\Pe\oplus\Pe^\vee\cong  {\F^\vee}^3\oplus \F^3 \cong\torus^3\oplus { \torus^\vee}^3$ as $\sts$-module.
\end{proof}

\begin{theorem} Let $X$ be a proper scheme over a perfect field $k$. Let $A_0$ be a central
simple division algebra over $k$.
  Let $k'/k$ be a separable cubic field extension of $k$ which splits $A_0$. Suppose it is Galois.
    For $\A=A_0\otimes \sts$, every left
 $\A$-module $\Pe$ of rank one with $N_\A(\Pe)\cong\sts$ is one of the following:\\
(i) $\Pe\cong \F\oplus\F\oplus \F$
for an indecomposable $\mathcal{O}_{X}$-module $\F$ of rank 3 with ${\rm det}\,\F\cong \mathcal{O}_{X}$,
\\ (ii) $\Pe\cong tr_{k'/k}(\F')$
for a suitable indecomposable $\mathcal{O}_{X'}$-module $\F'$ of rank 3 which is not defined over $X$
 with ${\rm det}\,\F'\cong \mathcal{O}_{X'}$,
\\ (iii) $\Pe\cong \lb\oplus\lb\oplus \lb\oplus \F\oplus\F\oplus \F$
for a line bundle $\lb $ and an indecomposable $\mathcal{O}_{X}$-module $\F$ of rank 2 with ${\rm det}\,\F\cong \lb^\vee$,
\\ (iv) $\Pe\cong tr_{k'/k}(\lb')\oplus \F\oplus\F\oplus \F$
for a line bundle $\lb' $ over $X'$ not defined over $X$ and an indecomposable  $\mathcal{O}_{X}$-module $\F$ of rank 2
with ${\rm det}(\F\otimes \mathcal{O}_{X'})\cong \lb'^\vee$,
\\ (v) $\Pe\cong \lb\oplus\lb\oplus \lb \oplus tr_{k'/k}(\F')$
for a line bundle $\lb $ over $X$  and an $\mathcal{O}_{X'}$-module $\F'$ of
 rank 2, which is not defined over $X$, with ${\rm det}\,\F'\cong (\lb^\vee)\otimes \mathcal{O}_{X'}$,
 \\ (vi) $\Pe\cong tr_{k'/k}(\lb')\oplus tr_{k'/k}(\F')$
for a line bundle $\lb' $ over $X'$  and an $\mathcal{O}_{X'}$-module $\F'$ of
 rank 2, both not defined over $X$, with ${\rm det}\,\F'\cong \lb'^\vee$,
 \\ (vii) $\Pe\cong \lb\oplus\lb\oplus \lb \oplus \s\oplus\s\oplus \s$
for two line bundles $\lb $, $\s$ over $X$ with $\lb\otimes \lb\otimes \lb\otimes
\s\otimes \s\otimes \s\cong  \mathcal{O}_{X}$,
 \\ (viii) $\Pe\cong \lb\oplus\lb\oplus \lb \oplus tr_{k'/k}(\s)$
for  line bundles $\lb $ over $X$, $\s'$ over $X'$, with $\s'$ not defined over $X$.
Again, $[(\lb \otimes \lb\otimes \lb)\otimes \mathcal{O}_{X'}]\otimes \s'\otimes \s'\otimes \s'\cong  \mathcal{O}_{X'}$,
 \\ (ix) $\Pe\cong tr_{k'/k}(\lb')\oplus tr_{k'/k}(\s')$
 for two line bundles $\lb' $, $\s'$ over $X'$ which are not defined over $X$, with $\lb'\otimes \lb'\otimes \lb'\otimes
\s'\otimes \s'\otimes \s'\cong  \mathcal{O}_{X'}$,
 \end{theorem}

\begin{proof} Since $\A$ is defined over $k$,
 there is a cubic form with adjoint and base point such that $\A^+=\J(N_\A,\sharp_\A,1)$.
Put $X'=X\times_k k'$ and let $\J=J(\A,\Pe,N)$. Then
$$\J\otimes_{\sts} \mathcal{O}_{X'}=J(\A\otimes \mathcal{O}_{X'},\Pe\otimes\mathcal{O}_{X'},N\otimes\mathcal{O}_{X'})$$
and $\A\otimes \mathcal{O}_{X'}\cong {\rm Mat}_3(\mathcal{O}_{X'})$.
Therefore there exists an $\mathcal{O}_{X'}$-module $\F'$ of rank 3 with ${\rm det}\,\F'\cong \mathcal{O}_{X'}$
such that
$$\Pe\otimes\mathcal{O}_{X'}\cong \F'\oplus\F'\oplus \F'.$$
This restricts the $\sts$-module structure of $\Pe$ and yields the assertion using [AEJ1] and a case-by-case study.

\end{proof}

So even in this very general scenario we can already exclude certain module structures for the left $\A$-modules of rank one
which carry a nondegenerate multiplicative form, if $\A=A_0\otimes \sts$ is as above. (If $\A=A_0\otimes \sts$ does not
contain a separable cubic extension which is Galois, a similar argument works using a Galois extension containing $k'$.)
This automatically restricts the module structures appearing in the first Tits construction starting with $\A$.
If one starts investigating special types of varieties, some of these cases will usually not occur, see Proposition .

\section{Brauer-Severi varieties}

 Let $k$ be a perfect field of characteristic not 2 or 3 and $X$ a Brauer-Severi variety with associated central simple
  algebra $D_0$ over $k$ of degree $3$. The Theorem of Krull-Schmidt for vector bundles
 holds over $X$ as well as its version for nondegenerate $\sts$-symmetric bilinear spaces ([AEJ1, 2] or
 [S, p.~274-275]).
 If $D_0$ is a division algebra and $k'/k$  a finite field extension which is a maximal subfield
 of $D_0$, then for $X'=X\times_k k'$ we have $X'\cong\mathbb{P}^2_{k'}$. In that case let $\I=tr_{k'/k}( \mathcal{O}_{X'}(-1))$
 be the canonical vector bundle of rank 3.  $\bigwedge^3\I\cong \omega_X$ is the canonical sheaf of $X$
 and  $\I\otimes \mathcal{O}_{X'}\cong  \mathcal{O}_{X'}(-1)^3$. There is an element $\lb$ generating the Picard group
 ${\rm Pic}\,X\cong\mathbb{Z}$ which satisfies $\lb\otimes  \mathcal{O}_{X'}\cong  \mathcal{O}_{X'}(3)$.
 Define $\lb(0)=\sts$, $\lb(m)=\lb\otimes\dots\otimes\lb$ ($m$-times) for $m>0$ and
 $\lb(m)=\lb^\vee\otimes\dots\otimes\lb^\vee$ ($(-m)$-times) for $m<0$ [Pu2]. In particular,
  $\omega_X\cong\lb(-1)$.
   Moreover,  $tr_{k'/k}( \mathcal{O}_{X'}(m))\otimes \mathcal{O}_{X'}\cong
   \mathcal{O}_{X'}(m)\oplus \mathcal{O}_{X'}(m)\oplus\mathcal{O}_{X'}(m)$
for all $m\in\mathbb{Z}$ and $^\sigma \mathcal{O}_{X'}(m)\cong \mathcal{O}_{X'}(m) $ for all $k'/k$-conjugates
[AEJ1, 1.9, 3.13].

For an algebraic field extension $l/k$, define $X_l=X\times_k l$.
We will repeatedly use that the canonical morphism
$${\rm Pic}\,\mathbb{P}_k^2\to {\rm Pic}\,\mathbb{P}_{l}^2,\,
\sts(m)\to \sts(m)\otimes \mathcal{O}_{X_l}$$
is bijective and that, for every algebraic field extension $l/k$ which does not split $X$, the canonical morphism
$${\rm Pic}\,X\to {\rm Pic}\,X_l,\,\lb(m)\to \lb(m)\otimes \mathcal{O}_{X_l}$$
is bijective.

  \begin{remark} A case-by-case study as in [Pu1, 3.3] shows that if $\mathcal{D}$ is a quaternion algebra over $X$ which
is neither split nor defined over $k$, then $\mathcal{D}\cong{\rm Quat} (\mathcal{F}, N)$ where $\mathcal{F}$ is
an absolutely indecomposable vector bundle on $X$ of rank 3 with trivial
determinant which carries a nondegenerate quadratic form. For example, choose $\mathcal{D}$ to be the quaternion algebra which is the extension  of $A = End_{{\Bbb H} [x, y]} (P)$
to ${\Bbb P}_{{\Bbb R}}^2$ with $P$ a non-free
projective ideal of ${\Bbb H} [x, y]$. Here, $x = \frac{Y}{Z}, \; y =
\frac{Y}{Z}$ where $(X, Y, Z)$ denote the homogeneous coordinates of ${\Bbb P}_{{\Bbb R}}^2$
(Knus, Parimala and Sridharan [K-P-S, 7.3]).
  \end{remark}

\begin{lemma} (i) Every first Tits construction starting
with $\mathcal{O}_X$ is defined over $k$.\\
(ii) Let $k_0/k$ be a separable  quadratic field extension with ${\rm Gal}(k_0/k)=\langle \omega\rangle$.
 Let $X_0=X\times_kk_0$. Every  Tits process over $X$ starting with $\B=\mathcal{O}_{X_0}$ and $*_\B=\omega$
 is defined over $k$.
\end{lemma}

\begin{proof} (i) Let $\J$ be a first Tits construction starting with $\mathcal{O}_X$.
If $X=\mathbb{P}_k^2$, every first Tits construction $\J$ over $\mathbb{P}_k^2$ starting with $\mathcal{O}_X$
is defined over $k$, since ${\rm Pic}(X)\cong\mathbb{Z}$.
 If $X$ is associated to a division algebra $D_0$ and if $k'$ is a splitting field of $D_0$, $X'=X\times_kk'$,
$\J\otimes_k k'$ is a first Tits construction over $X'=\mathbb{P}_k^2$ and as such globally free as
$\mathcal{O}_{X'}$-module. Thus so is $\J$, and by Lemma 2 it must be defined over $k$.
\\ (ii) Let $X=\mathbb{P}_k^2$. If $\Pe\in{\rm Pic}\,X_0$ is $\sts$-admissible then $\Pe=\mathcal{O}_{X_0}$
 (each $\Pe$ is already defined over $X$, so $^\omega\Pe\cong\Pe$ and thus $^\omega\Pe\cong\Pe^\vee$ implies
 $\Pe\cong\Pe^\vee$). Therefore
 $\J=\J(\mathcal{O}_{X_0},\sts,\Pe,N,*)$ is globally free as $\sts$-module implying that it is defined over $k$
 by Lemma 2.
 If the central simple algebra $D_0$ associated to $X$ is a division algebra and $\J$ is a  Tits process over
  $X$, then, given any splitting field $l$ of $D_0$, $\J\otimes \mathcal{O}_{X_l}$ is  a Tits process over
$\mathcal{O}_{X_l}$ and thus defined over $l$. Hence $\J$ is globally free as $\sts$-module
implying that it is defined over $k$.
\end{proof}

\subsection{ \'Etale first Tits constructions}

\begin{theorem} Let $k'/k$ be a  cubic field extension of $k$.
Every  \'etale first Tits constructions starting with $\A=k'\otimes \sts$ is defined over $k$.
\end{theorem}

\begin{proof}  Define $X'=X\times_k k'$ to be the Brauer-Severi variety obtained from $X$ by change of bases.
Identify ${\rm Pic}\,\A$ and $ {\rm Pic}\,X'$.
Suppose first that $X$ is rational. Then the map
 $$ {\rm Pic}\,X\to {\rm Pic}\,\A,\sts(m)\to \sts(m)\otimes \mathcal{O}_{X'}$$
 is bijective. Therefore every $\Pe\in {\rm Pic}\,\A$ is of the form $\Pe= \mathcal{O}_{X'}(m)
\cong \sts(m)\otimes \mathcal{O}_{X'}\cong (k'\otimes \sts)\otimes \sts(m)$.
 If $\Pe\cong (k'\otimes \sts)\otimes \sts(m)$
 for some $m\in\mathbb{Z}$ and $N$ is the norm of $\A$, then
$$\sts\cong N(\Pe)\cong N((l\otimes \sts)\otimes \sts(m))\cong\sts(3m)$$
implies $m=0$.
Thus there exists only the classical first Tits construction for $\A$ which yields a Jordan algebra of rank 9
which is defined over $k$.

If $X$ is associated to a division algebra $D_0$ over $k$, we distinguish two cases:
Suppose first that $k'$ is no splitting field of $D_0$.
Then the map $$ {\rm Pic}\,X\to {\rm Pic}\,\A,\lb(m)\to \lb(m)\otimes \mathcal{O}_{X'}$$
 is bijective and an analogous argument as in the previous case shows that
 there exists only the classical first Tits construction for $\A$, which yields a Jordan algebra
which is defined over $k$.\\
Suppose now that $k'$ is a splitting field of $D_0$. It suffices to consider the case that $k'$ is a maximal subfield
of $D_0$.

 Let $k^+$ be a Galois extension of $k$ containing $k'$ with ${\rm Gal}(k^+/k)=< id,\sigma_1,\sigma_2 >$. Since
 $X'=\mathbb{P}_{k'}^2$, there are line bundles $\mathcal{O}_{X'}(3m\pm 1)$ over $X'$ which are not
 defined over $X$ ($m\in\mathbb{Z}$). However,
$$\mathcal{O}_{X^+}(3m\pm 1)^+\otimes\, ^{\sigma_1}\mathcal{O}_{X^+}(3m\pm 1)^+\otimes \, ^{\sigma_2}\mathcal{O}_{X^+}
(3m\pm 1)^+\cong  \mathcal{O}_{X^+}(3(3m\pm 1)) \cong \mathcal{O}_{X^+}$$
yields the contradiction that there exists $m\in \mathbb{Z}$ such that $3m\pm 1=0$, so the case
$$\Pe\cong tr_{k'/k}(\mathcal{O}_{X'}(3m\pm 1))$$
cannot occur.
\end{proof}

\begin{lemma} Let $\torus$ be a non-split quadratic  \'etale algebra over $X$.  Let $\A=\sts\times\torus$.
Then $\J(\A,\Pe,N)$ is a Freudenthal algebra over $X$ of rank 9 with underlying module structure\\
 $$\J=\sts^3\oplus [ \sts(-2m)\oplus \sts(m)\oplus\sts(m)]\oplus [ \sts(2m)\oplus \sts(-m)\oplus\sts(-m)] $$
 if $X=\mathbb{P}_k^2$ and
 $$\J=\sts^3\oplus [ \lb(-2m)\oplus \lb(m)\oplus\lb(m)]\oplus  [\lb(2m)\oplus \lb(-m)\oplus\lb(-m)]$$
  if $X$ is associated to a division algebra $D_0$.
\end{lemma}

\begin{proof}
We have  $\torus=k'\otimes \sts$ with $k'$ a quadratic field extension of $k$ and ${\rm Pic}\,\torus
\cong\mathbb{Z}$.  Let ${\rm Gal}(k'/k)=\langle \sigma\rangle$.
By Lemma 5, every left $\A$-module $\Pe$ of rank one with $N_\A(\Pe)\cong\sts$ is
  a direct sum of line bundles $\lb_0,\M_0,\s_0$ over $X$ with
$\lb_0\otimes\M_0\otimes\s_0\cong\mathcal{O}_{X}$; i.e.,
 $$\Pe\cong \lb_0\oplus\M_0\oplus \s_0,$$
  or there is a line bundle $\s_0$ over $X$ and a line bundle $\lb$  over $X'$ which is not defined over $X$
satisfying $\lb\otimes\, ^\sigma\lb\otimes\s_0\cong\mathcal{O}_{X'}$
and
$$\Pe\cong tr_{k'/k}(\lb)\oplus\s_0.$$
As we will see, this last case cannot  happen:

Since $\check{H}^1(X,\A^\times)={\rm Pic}\,X\oplus {\rm Pic}\,\torus$, every left $\A$-module $\Pe$
of rank one satisfies $\Pe\cong \lb\oplus\F$ with $\F\in {\rm Pic}\torus$ and an invertible $\sts$-module $\lb$.
 We have $ N_\A(\Pe)\cong \sts$ iff $\lb\otimes N_\torus(\F)\cong \sts$ iff $\lb^\vee\cong N_\torus(\F)$.
For $\M\in {\rm Pic}\,X$, $\M\otimes\torus\in {\rm Pic}\,\torus$ has $N_\torus(\M\otimes\torus)\cong \M
^2$.

Now ${\rm Pic}\torus$ can  be identified with  ${\rm Pic}X'$ and  the canonical morphism
 $$G: {\rm Pic}X'\to {\rm Pic}X', \lb\to \lb\otimes \mathcal{O}_{X'}$$
 is bijective.
\\{\it Case I.} Suppose $X=\mathbb{P}_k^2$, then  each $\Pe_0\in{\rm Pic}\torus$ is of the form
  $\Pe_0\cong (k'\otimes \sts)\otimes \sts(m)$. If $N$ is the norm on $k'\otimes_k \sts$ then
$$ N(\Pe_0)\cong N((k'\otimes \sts)\otimes \sts(m))\cong\sts(2m)$$
Therefore,
$$\Pe\cong \sts(-2m)\oplus \sts(m)\oplus\sts(m)$$
as $\sts$-module for some $m\in\mathbb{Z}$.
\\ {\it Case II.} Suppose $X$ is associated to a division algebra $D_0$, then each $\Pe_0\in{\rm Pic}\torus$ is of the form
  $\Pe_0\cong (k'\otimes \sts)\otimes \lb(m)$.  If $N$ is the norm on $k'\otimes_k \sts$ then
$$ N(\Pe_0)\cong N((k'\otimes \sts)\otimes \lb(m))\cong\lb(2m).$$
Hence
$$\Pe\cong \lb(-2m)\oplus \lb(m)\oplus\lb(m)$$
as $\sts$-module for some $m\in\mathbb{Z}$.
\end{proof}

\subsection{Albert algebras over $X$ obtained by first Tits constructions}

  Suppose $D_0$ is a division algebra.

\begin{example}
Starting with the Azumaya algebra $\A=\E nd_X(\I)\cong D^{\rm op}\otimes_k \sts$ defined using
the opposite algebra of $D_0$, we compute some possible first Tits
 constructions (the indecomposable vector bundle $\E=\lb(s)\otimes \I$  yields the same Azumaya
 algebra $\E nd_X(\E)=(\lb(s)\otimes \I)\otimes
(\lb(-s)\otimes \I^\vee)\cong \E nd_X(\I)$): All locally free left $\A$-modules of rank 1 admitting a multiplicative norm
  must be of the kind $\h om_X(\F,\I)$ with $\F$ a vector bundle of rank 3 such that $\bigwedge^3\F\cong \lb(-1)$.
\\ (1) If $\F$ is a direct sum of line bundles then
 $$\F=\lb(m_1)\oplus\lb(m_2)\oplus\lb(-m_1-m_2-1)$$
 with arbitrary $m_1,m_2\in\mathbb{Z}$.
 $\J(D_0^{\rm op}\otimes \sts,\Pe,N)$ with $\Pe=\h om_X(\F,\I)$ is an Albert algebra over $X$ isomorphic to

\begin{align*}
\mathcal{O}_X^9
 &\oplus
\left [\begin {array}{ccc}
 tr_{k'/k}( \mathcal{O}_{X'}(-1-3m_1)) &  tr_{k'/k}( \mathcal{O}_{X'}(-1-3m_2))  &  tr_{k'/k}( \mathcal{O}_{X'}
 (2+3m_1+3m_2))\\
\end {array}\right ]
 \\&
\oplus
\left [\begin {array}{ccc}
 tr_{k'/k}( \mathcal{O}_{X'}(1+3m_1)) \\
 tr_{k'/k}( \mathcal{O}_{X'}(1+3m_2))  \\
 tr_{k'/k}( \mathcal{O}_{X'}(-2-3m_1-3m_2)) \\
\end {array}\right ]
\end{align*}
as $\sts$-module.
\\(2) Every indecomposable (but not absolutely indecomposable) $\F$ of rank 3 is isomorphic to
$\lb(m)\otimes tr_{k'/k}( \mathcal{O}_{X'}(\pm 1))\cong tr_{k'/k}( \mathcal{O}_{X'}(3m \pm 1))$
  and has determinant
  $\bigwedge^3\F\cong \lb(m\pm 1)$, since
 $$\bigwedge\,^3(\F\otimes  \mathcal{O}_{X'})\cong \mathcal{O}_{X'}(-1+3m)\otimes \mathcal{O}_{X'}(-1+3m)\otimes
 \mathcal{O}_{X'}(3m \pm 1)=\mathcal{O}_{X'}(3(3m\pm 1)).$$
  Hence the only one with suitable determinant is $\I$
 itself which yields the classical first Tits construction which is defined over $k$.\\
(3) It is not possible that $\F$ decomposes as a direct sum of an indecomposable (but not absolutely indecomposable)
 bundle and a line bundle: If this were true then we could write $\F=tr_{l_k}(\M)\oplus \dots$ with a line bundle
 $\M$ over $X\otimes_k l$, where $l$ is a quadratic field extension of $k$. This, however, would not yield an
 indecomposable bundle of rank 2 over $X$, since each line bundle $\M$ over $X_l$ is already defined over $X$,
 a contradiction.\\
(4) It may be possible that $\F$ is an absolutely indecomposable bundle of rank 3, or the direct sum of an absolutely
 indecomposable bundle of rank 2 and a line bundle. We know too little about the vector bundles over $X$
 to treat these cases.
\end{example}

\begin{example} Let $\A=\E nd_X(\E)$ with $\E=\lb(m_1)\oplus\lb(m_2)\oplus\lb(m_3)$.
\\(1) If $\F$ is a direct sum of line bundles then
 $$\F=\lb(l_1)\oplus\lb(l_2)\oplus\lb(l_3)$$
with $m_1+m_2+m_3=l_1+l_2+l_3$.
\\  (2) If $\F$ is indecomposable (but not absolutely so)
then we may choose
$$\F= tr_{k'/k}( \mathcal{O}_{X'}(3m - 1))$$
with $-1+m=m_1+m_2+m_3$, i.e., $m=m_1+m_2+m_3+1$.
 Put $\A=\E nd_X(\E)$, then $J(\A,\Pe,N)$, $\Pe=\h om_X(\F,\E)$ is an Albert algebra over $X$ and,
 for $\F$ a direct sum of line bundles, isomorphic to
\begin{align*}
\left [\begin {array}{ccc}
\sts & \lb(a) & \lb(b)\\
\lb(-a) &  \sts & \lb(b-a)  \\
\lb(-b) & \lb(a-b) &  \sts \\
\end {array}\right ]
& \oplus
\left [\begin {array}{ccc}
\lb(c) & \lb(a+d)  & \lb(b-c-d)\\
\lb(c-a) & \lb(d) & \lb(b-a-c-d)  \\
\lb(c-b)  &  \lb(a-b+d)  &  \lb(-c-d)  \\
\end {array}\right ]
\\ & \oplus
\left [\begin {array}{ccc}
\lb(-c) & \lb(a-c) &  \lb(b-c) \\
\lb(-a-d)  &  \lb(-d) &  \lb(b-a-d) \\
 \lb(-b+c+d)  &  \lb(a-b+c+d) &  \lb(c+d) \\
\end {array}\right ]
\end{align*}
with $a=m_1-m_2$, $b=m_1-m_3$, $c=m_1-l_1$ and $d=m_2-l_2$. This algebra is of the same type as the one computed
 in [Pu3, Example 3] for $X=\mathbb{P}_k^1$. Indeed, Albert algebras of this type
can be constructed over $\mathbb{P}_k^n$ for any integer $n$.
\\(3) If $\F= tr_{k'/k}( \mathcal{O}_{X'}(3m-1))$ then $\A$ is isomorphic to
\begin{align*}
& \left [\begin {array}{ccc}
\sts & \lb(a) & \lb(b)\\
\lb(-a) &  \sts & \lb(b-a)  \\
\lb(-b) & \lb(a-b) &  \sts \\
\end {array}\right ]
  \oplus
\left [\begin {array}{ccc}
tr_{k'/k}( \mathcal{O}_{X'}(-2-3m_2-3m_3)) \\
tr_{k'/k}( \mathcal{O}_{X'}(-2-3m_1-3m_3)) \\
tr_{k'/k}( \mathcal{O}_{X'}(-2-3m_1-3m_2))  \\
\end {array}\right ]
\\ &\oplus
\left [\begin {array}{ccc}
tr_{k'/k}( \mathcal{O}_{X'}(2+3m_2+3m_3)) & \mathcal{O}_{X'}(2+3m_1+3m_3)) &  \mathcal{O}_{X'}(2+3m_1+3m_2)) \\
\end {array}\right ]
\end{align*}
as $\sts$-module.
\end{example}

\begin{proposition} Let $A_0$ be a central simple division algebra over $k$.
Let $k'/k$ be a  cubic field extension of $k$ which splits $A_0$. Suppose it is Galois.
For $\A=A_0\otimes \sts$, every left $\A$-module $\Pe$ of rank one with $N_\A(\Pe)\cong\sts$ is of one of the types described
in Theorem 5. If $X=\mathbb{P}_k^2$ or if $D_0$ is a division algebra over $k$ and $A_0$ and $D_0$ have a common splitting
field $k'$ as above, then the cases (iv), (vi), (viii) and (ix) do not occur.
\end{proposition}

\begin{proof} This follows from Theorem 5.
\end{proof}

\begin{example}
Let us look at examples of Albert algebras of the kind $\h_3 (\mathcal{C},\Gamma)$ with $\mathcal{C}=
{\rm Zor}\,(\torus,\alpha)$ a split octonion algebra over $X$.
For trivial $\Gamma$, $\mathcal{H}_3(\comp,\Gamma)$ is a
 first Tits construction starting with ${\rm Mat}_3(\sts)$.\\
Let $\torus$ be a vector bundle of constant rank 3 over $X$ such that $\bigwedge^3\torus\cong\sts$.
\\ (1) If $\torus$ is the direct sum of line bundles, then
$$\torus=\lb(m_1)\oplus\lb(m_2)\oplus\lb(-(m_1+m_2))$$
with $m_1\geq m_2\geq 0$. We obtain
\[{\rm Zor}\,(\torus,\alpha)=
\left [\begin {array}{ccc}
\sts & \lb(m_1)\oplus\lb(m_2)\oplus\lb(-(m_1+m_2)) \\
 \lb(-m_1)\oplus\lb(-m_2)\oplus\lb(m_1+m_2) & \sts  \\
\end {array}
\right ].
\]
Hence $\mathcal{H}_3(\comp,\Gamma)$ has
$$ \sts^9\oplus \lb(m_1)\oplus[\lb(m_2)\oplus\lb(-(m_1+m_2))]^3 \oplus [  \lb(-m_1)\oplus\lb(-m_2)\oplus\lb(m_1+m_2)]^3 $$
as underlying module structure.
\\ (2) The only indecomposable (but not absolutely indecomposable) $\torus$ of rank 3 with trivial determinant
 is isomorphic to a
vector bundle of the form
$$\lb\otimes tr_{k'/k}( \mathcal{O}_{X'}(-1))\cong tr_{k'/k}( \mathcal{O}_{X'}(3m-1))$$
or to its dual, and yields
\[{\rm Zor}\,(\torus,\alpha)=
\left [\begin {array}{ccc}
\sts & tr_{k'/k}( \mathcal{O}_{X'}(-1+3m)) \\
 tr_{k'/k}( \mathcal{O}_{X'}(1-3m)) & \sts  \\
\end {array}
\right ]
\]
Hence $\mathcal{H}_3(\comp,\Gamma)$ has
$$\sts^9\oplus  tr_{k'/k}( \mathcal{O}_{X'}(-1+3m))^3 \oplus tr_{k'/k}( \mathcal{O}_{X'}(1-3m))^3 $$
as underlying module structure.
\\ (3) Over ${\Bbb P}_{{\Bbb R}}^2$, there exist absolutely indecomposable vector bundles  $\mathcal{T}$ on $X$ of rank 3
 with trivial determinant.

  For example, choose $\torus$ to be the trace zero elements of the quaternion algebra $\mathcal{D}$
  in Remark 9.
  Then $\mathcal{H}_3(\comp,\Gamma)\cong \sts^9\oplus \torus^3\oplus  \torus^{\vee 3}$ as $\sts$-module.

\end{example}

\subsection{Albert algebras over $X$ obtained by the Tits process}

We finish by looking at some examples of Azumaya algebras which can be used for a Tits process:

\begin{theorem}  Let $D_0$ be a division algebra  over  $k$.
 Let $k'$ be a separable quadratic field extension of $k$ with ${\rm Gal}(k'/k)=\langle\omega\rangle$.
  Let $X'=X\times_kk'$.  Then
$$\B=\E nd_{X'} (\mathcal{O}_{X'}\oplus\lb(m)\oplus\lb(-mP))$$
 is an Azumaya algebra of rank 9 over $X'$ permitting an involution $*_\B$ adjoint to the hermitian form defined on
 $\mathcal{O}_{X'}\oplus\lb(m)\oplus\lb(-m)$, such that $\h(\mathcal{O}_{X'},*_\B)=\sts$
and $\B$, $\mathcal{O}_X'=\mathcal{O}_{X'}$ and $*_\B$
 are suitable for the Tits process $\J(\B, \h (\B,*_\B),\Pe,N,*)$.
 We have
$$ \h (\B,*_\B)\cong \mathcal{O}_{X}^3\oplus \lb(m)^2\oplus \lb(-m)^2\oplus \lb(2m)\oplus \lb(-2m)$$
as underlying $\sts$-module structure.
The $\sts$-module structure of $\Pe$ is isomorphic to one of the following:\\
(i) $$\begin{array}{l}
\Pe\cong \lb(n_1)\oplus\lb(n_2)\oplus\lb(-(n_1+n_2))\\
\oplus\lb(-m+n_1)\oplus\lb((-m+n_2))\oplus\lb(-m-(n_1+n_2))\\
\oplus\lb(m+n_1)\oplus\lb(m+n_2)\oplus\lb(m-(n_1+n_2))\\
\oplus
\lb(-n_1)\oplus\lb(-n_2)\oplus\lb(n_1+n_2)\\
\oplus\lb((m-n_1))\oplus\lb((m-n_2))\oplus\lb(m+n_1+n_2)\\
\oplus\lb(-m-n_1)\oplus\lb(-m-n_2)\oplus\lb(-m+n_1+n_2).
\end{array}$$
for arbitrarily chosen $n_i\in\mathbb{Z}$.
\\ (ii) \begin{align*}\Pe\cong  \lb(-n) \oplus tr_{k'/k}(\M^\vee) \oplus\lb(m-n) \oplus
 \lb(m)\otimes tr_{k'/k}(\M^\vee)
\oplus \lb(-m-n) \oplus\\
 \lb(-m)\otimes tr_{k'/k}(\M^\vee) \oplus
 \lb(n) \oplus \lb(-m+n) \oplus \lb(m+n)
 \end{align*}
for an arbitrarily chosen $n\in\mathbb{Z}$
 and for an absolutely indecomposable vector bundle $\M$ on $X$ of rank 2 which is not defined over $X$ and satisfies
$^\omega\M\cong\M^\vee$ and ${\rm det}\,\M\cong\lb(-n)$.
\\(iii) \begin{align*} \Pe\cong [\lb(-n) \oplus \M_0^\vee \oplus\lb(m-n) \oplus \lb(m)\otimes \M_0^\vee
\oplus \lb(-m-n) \oplus \lb(-m)\otimes \M_0^\vee ]\oplus \\
[ \lb(n) \oplus \M_0 \oplus\lb(-m+n) \oplus \lb(-m)\otimes \M_0
\oplus \lb(m+n) \oplus \lb(m)\otimes \M_0]
\end{align*}
for an arbitrarily chosen $n\in\mathbb{Z}$
 and for an absolutely indecomposable vector bundle $\M_0$ on $X'$ of rank 2 such that ${\rm det}\,\M\cong\lb(-n)\otimes
\mathcal{O}_{X'}$.
\\ (iv)
\[\Pe\cong [ \G_0^\vee\oplus \lb(m)\otimes\G_0^\vee\oplus\lb(-m) \otimes\G_0^\vee ]\oplus
[ \G_0\oplus \lb(-m)\otimes\G_0\oplus\lb(m) \otimes\G_0] \]
for an absolutely indecomposable vector bundle $\G_0$ on $X'$ of rank 3 over $X$ with trivial determinant.
\\(v)
\[\Pe\cong  tr_{k'/k}( \G)\oplus \lb(m)\otimes  tr_{k'/k}( \G)\oplus\lb(-m)\otimes  tr_{k'/k}( \G) \]
for an absolutely indecomposable vector bundle $\G$ on $X'$ of rank 3 with trivial determinant
 which is not defined over $X$ and satisfies $^\omega\G\cong\G^\vee$.
 \\ (vi) Suppose there is a vector bundle $\G$ of rank 3 which is indecomposable, but not absolutely so, i.e.,
$\G\cong  \lb(n)\otimes tr_{l'/k'}(\mathcal{O}_{X_{l'}}(\pm 1))\cong tr_{l'/k'}(\mathcal{O}_{X_{l'}}(3n\pm1))$
for a separable cubic field extension $l'/k'$ which splits the central simple division algebra $D_0$ associated with $X'$
($n\in\mathbb{Z}$). If it has trivial determinant (we do not know if this can happen), then
 there is a cubic separable field extension $l/k$ which is a splitting field of $D_0$ such that
 \begin{align*}
 \Pe\cong [ tr_{l/k}(\mathcal{O}_{X_{l}}(3n\mp 1))\oplus tr_{l/k}(\mathcal{O}_{X_{l}}
(3(n+m)\mp 1))\oplus
tr_{l/k}(\mathcal{O}_{X_{l}}(3(n-m)\mp 1))]\\ \oplus
[ tr_{l/k}(\mathcal{O}_{X_{l}}(3n\pm 1))\oplus
 tr_{l/k}(\mathcal{O}_{X_{l}}(3(n-m)\pm 1))\oplus tr_{l/k}(\mathcal{O}_{X_{l}}(3(n+m)\pm 1))].
  \end{align*}
  or
 \begin{align*}
 \Pe\cong  tr_{l'/k}(\mathcal{O}_{X_{l}}(3n\mp 1))\oplus tr_{l'/k}(\mathcal{O}_{X_{l}}
(3(n+m)\mp 1))\oplus
tr_{l'/k}(\mathcal{O}_{X_{l}}(3(n-m)\mp 1)).
  \end{align*}
\end{theorem}

\begin{proof} $X'$ is nonrational as well. Since
$\J(\B, \h (\B,*_\B),\Pe,N,*)\otimes  \mathcal{O}_{X'} $ contains $\B^+$ it is a first Tits construction
$\J(\B,\Pe_0,N_0)$ where
$$\Pe_0\cong\h om_{\mathcal{O}_{X'}}(\G,\mathcal{O}_{X'}\oplus \lb(m)\oplus \lb(-m))$$
for a vector bundle $\G$ of rank 3 with ${\rm det}\,\G\cong \mathcal{O}_{X'}$. We distinguish the following cases:
\\ (1) $\G\cong \lb(n_1)\oplus \lb(n_2)\oplus \lb(-(n_1+n_2))$. This implies (i).
\\ (2) If $\G\cong  \lb(n)\oplus \M$ with $\M$ an indecomposable
 vector bundle of constant rank 2 then $\M$ must be absolutely indecomposable, since every line bundle over $X_{l'}$
 with $l'/k'$ a separable quadratic field extension is already defined over $X'$. Moreover, we need
 ${\rm det}\,\M\cong \lb(-m)$. We get

 \[ \Pe_0\cong  \h om_{\mathcal{O}_{X'}}( \lb(n)\oplus \M,\mathcal{O}_{X'}\oplus \lb(m)\oplus \lb(-m))\cong\\
\left [\begin {array}{cc}
\lb(-n) & \M^\vee \\
\lb(m-n) & \lb(m)\otimes \M^\vee \\
\lb(-m-n) & \lb(-m)\otimes \M^\vee\\
\end {array}
\right ].
\]

Thus
\begin{align*}\Pe\otimes\mathcal{O}_{X'}\cong [\lb(-n) \oplus \M^\vee \oplus\lb(m-n) \oplus \lb(m)\otimes \M^\vee
\oplus \lb(-m-n) \oplus \lb(-m)\otimes \M^\vee ] \\
\oplus
[ \lb(n) \oplus \M \oplus\lb(-m+n) \oplus \lb(-m)\otimes \M
\oplus \lb(m+n) \oplus \lb(m)\otimes \M]
 \end{align*}
If $\M$ is not already defined over $X$ (note that since $^\omega\M\cong\M^\vee$ by Proposition 3,
$ tr_{k'/k}( \M)\cong tr_{k'/k}( \M^\vee)$) then we get (ii), i.e.,
\begin{align*} \lb(-n) \oplus tr_{k'/k}(\M^\vee) \oplus\lb(m-n) \oplus
 \lb(m)\otimes tr_{k'/k}(\M^\vee)
\oplus \lb(-m-n) \oplus\\
 \lb(-m)\otimes tr_{k'/k}(\M^\vee) \oplus
 \lb(n) \oplus \lb(-m+n) \oplus \lb(m+n)
 \end{align*}
as a possible $\sts$-module structure of $\Pe$.\\
If $\M$ is already defined over $X$, that is $\M=\M_0\otimes
\mathcal{O}_{X'}$ then we get
\begin{align*} [\lb(-n) \oplus \M_0^\vee \oplus\lb(m-n) \oplus \lb(m)\otimes \M_0^\vee
\oplus \lb(-m-n) \oplus \lb(-m)\otimes \M_0^\vee ]\oplus \\
[ \lb(n) \oplus \M_0 \oplus\lb(-m+n) \oplus \lb(-m)\otimes \M_0
\oplus \lb(m+n) \oplus \lb(m)\otimes \M_0]
\end{align*}
as a possible $\sts$-module structure of $\Pe$, i.e., (iii).
 \\ (3) If $\G$ is absolutely indecomposable with trivial determinant (this case exists!), then
  \[ \Pe_0\cong  \h om_{\mathcal{O}_{X'}}( \G,\mathcal{O}_{X'}\oplus \lb(m)\oplus \lb(-m))\cong\\
    \G^\vee\oplus \lb(m)\otimes\G^\vee\oplus\lb(-m) \otimes\G^\vee \]
    and
\[\Pe\otimes\mathcal{O}_{X'}\cong [ \G^\vee\oplus \lb(m)\otimes\G^\vee\oplus\lb(-m) \otimes\G^\vee ]\oplus
[ \G\oplus \lb(-m)\otimes\G\oplus\lb(m) \otimes\G ]. \]
If $\G$ is already defined over $X$, that is $\G=\G_0\otimes
\mathcal{O}_{X'}$ then we get (iv):

\[\Pe\cong [ \G_0^\vee\oplus \lb(m)\otimes\G_0^\vee\oplus\lb(-m) \otimes\G_0^\vee ]\oplus
[ \G_0\oplus \lb(-m)\otimes\G_0\oplus\lb(m) \otimes\G_0]. \]
If $\G$ is not already defined over $X$ then (since $^\omega\G\cong\G^\vee$ by Proposition 3,
$ tr_{k'/k}( \G)\cong tr_{k'/k}( \G^\vee)$) we get (v):

\[ tr_{k'/k}( \G)\oplus \lb(m)\otimes  tr_{k'/k}( \G)\oplus\lb(-m)\otimes  tr_{k'/k}( \G). \]
\\(4) $\G$ is indecomposable, but not absolutely so, i.e.,
$\G\cong  \lb(n)\otimes tr_{l'/k'}(\mathcal{O}_{X_{l'}}(\pm 1))\cong tr_{l'/k'}(\mathcal{O}_{X_{l'}}(3n\pm1))$
for a separable cubic field extension $l'/k'$ which splits the central simple division algebra associated with $X'$
($n_i\in\mathbb{Z}$).

We have $\G\otimes \mathcal{O}_{X_{l'}}\cong  \mathcal{O}_{X_{l'}}(3n\pm 1)\oplus  ^{\sigma_1}
\mathcal{O}_{X_{l'}}(3n\pm 1) \oplus
^{\sigma_2}\mathcal{O}_{X_{l'}}(3n\pm 1)$. If $\G\otimes \mathcal{O}_{X_{l'}}$ has trivial determinant
(and we do not know if this can happen) then
$$\Pe_0
\cong tr_{l'/k'}(\mathcal{O}_{X_{l'}}(3n\mp 1))\oplus tr_{l'/k'}(\mathcal{O}_{X_{l'}}(3(n+m)\mp 1))\oplus
tr_{l'/k'}(\mathcal{O}_{X_{l'}}(3(n-m)\mp 1)) $$
and so
 \begin{align*}
 \Pe\otimes \mathcal{O}_{X'}\cong [ tr_{l'/k'}(\mathcal{O}_{X_{l'}}(3n\mp 1))\oplus tr_{l'/k'}(\mathcal{O}_{X_{l'}}
(3(n+m)\mp 1))\oplus
tr_{l'/k'}(\mathcal{O}_{X_{l'}}(3(n-m)\mp 1))]\\ \oplus
[ tr_{l'/k'}(\mathcal{O}_{X_{l'}}(3n\pm 1))\oplus
 tr_{l'/k'}(\mathcal{O}_{X_{l'}}(3(n-m)\pm 1))\oplus tr_{l'/k'}(\mathcal{O}_{X_{l'}}(3(n+m)\pm 1))].
  \end{align*}
These imply the assertion using [AEJ1] and the fact that every line bundle over $X'$ is already defined over $X$ here.
\end{proof}

\begin{theorem}  Let $X=\mathbb{P}_k^2$.
 Let $k'$ be a separable quadratic field extension of $k$ with ${\rm Gal}(k'/k)=\langle\omega\rangle$.
  Let $X'=X\times_kk'$.  Then
$$\B=\E nd_{X'} (\mathcal{O}_{X'}\oplus\mathcal{O}_{X'}(m)\oplus\mathcal{O}_{X'}(-m))$$
 is an Azumaya algebra of rank 9 over $X'$ permitting an involution $*_\B$ adjoint to the hermitian form defined on
 $\mathcal{O}_{X'}\oplus\mathcal{O}_{X'}(m)\oplus\mathcal{O}_{X'}(-m)$ such that $\h(\mathcal{O}_{X'},*_\B)=\sts$
and $\B$, $\mathcal{O}_X'=\mathcal{O}_{X'}$ and $*_\B$
 are suitable for the Tits process $\J(\B, \h (\B,*_\B),\Pe,N,*)$.
 We have
$$ \h (\B,*_\B)\cong \mathcal{O}_{X}^3\oplus \mathcal{O}_{X}(m)^2\oplus \mathcal{O}_{X}(-m)^2\oplus \mathcal{O}_{X}(2m)
\oplus \mathcal{O}_{X}(-2m)$$
as underlying $\sts$-module structure.
The $\sts$-module structure of $\Pe$ is isomorphic to one of the following:\\
(i) $$\begin{array}{l}
\Pe\cong\sts(n_1)\oplus\sts(n_2)\oplus\sts(-(n_1+n_2))\\
\oplus\sts(-m+n_1)\oplus\sts((-m+n_2))\oplus\sts(-m-(n_1+n_2))\\
\oplus\sts(m+n_1)\oplus\sts(m+n_2)\oplus\sts(m-(n_1+n_2))\\
\oplus
\sts(-n_1)\oplus\sts(-n_2)\oplus\sts(n_1+n_2)\\
\oplus\sts((m-n_1))\oplus\sts((m-n_2))\oplus\sts(m+n_1+n_2)\\
\oplus\sts(-m-n_1)\oplus\sts(-m-n_2)\oplus\sts(-m+n_1+n_2).
\end{array}$$
for arbitrarily chosen $n_i\in\mathbb{Z}$.
\\ (ii) \begin{align*}\Pe\cong\sts(-n) \oplus tr_{k'/k}(\M^\vee) \oplus\sts(m-n) \oplus
\sts(m)\otimes tr_{k'/k}(\M^\vee)
\oplus\sts(-m-n) \oplus\\
\sts(-m)\otimes tr_{k'/k}(\M^\vee) \oplus
\sts(n) \oplus\sts(-m+n) \oplus\sts(m+n)
 \end{align*}
for an arbitrarily chosen $n\in\mathbb{Z}$ and for an absolutely indecomposable vector bundle $\M$ on $X'$
 of rank 2 which is not defined over $X$ and satisfies
$^\omega\M\cong\M^\vee$ and ${\rm det}\,\M\cong\mathcal{O}_{X'}(-n)$.
\\(iii) \begin{align*} \Pe\cong[\sts(-n) \oplus \M_0^\vee \oplus\sts(m-n) \oplus\sts(m)\otimes \M_0^\vee
\oplus\sts(-m-n) \oplus\sts(-m)\otimes \M_0^\vee ]\oplus \\
[\sts(n) \oplus \M_0 \oplus\sts(-m+n) \oplus\sts(-m)\otimes \M_0
\oplus\sts(m+n) \oplus\sts(m)\otimes \M_0]
\end{align*}
for an arbitrarily chosen $n\in\mathbb{Z}$ and for an absolutely indecomposable vector bundle
$\M_0$ on $X'$ of rank 2 such that ${\rm det}\,\M\cong\sts(-n)\otimes
\mathcal{O}_{X'}$.
\\ (iv)
\[\Pe\cong[ \G_0^\vee\oplus\sts(m)\otimes\G_0^\vee\oplus\sts(-m) \otimes\G_0^\vee ]\oplus
[ \G_0\oplus\sts(-m)\otimes\G_0\oplus\sts(m) \otimes\G_0] \]
for an absolutely indecomposable vector bundle $\G_0$ on $X'$ of rank 3 over $X$ with trivial determinant.
\\(v)
\[ \Pe\cong tr_{k'/k}( \G)\oplus\sts(m)\otimes  tr_{k'/k}( \G)\oplus\sts(-m)\otimes  tr_{k'/k}( \G) \]
for an absolutely indecomposable vector bundle $\G$ on $X'$ of rank 3 with trivial determinant
 which is not defined over $X$ and satisfies $^\omega\G\cong\G^\vee$.
\end{theorem}

\begin{proof}
$\J(\B, \h (\B,*_\B),\Pe,N,*)\otimes  \mathcal{O}_{X'} $ contains $\B^+$. Hence it is a first Tits construction
$\J(\B,\Pe_0,N_0)$, where
$$\Pe_0\cong\h om_{\mathcal{O}_{X'}}(\G,\mathcal{O}_{X'}\oplus \mathcal{O}_{X'}(m)\oplus \mathcal{O}_{X'}(-m))$$
for a vector bundle $\G$ of rank 3 with ${\rm det}\,\G\cong \mathcal{O}_{X'}$. We distinguish the following cases:
\\ (1) $\G\cong \mathcal{O}_{X'}(n_1)\oplus \mathcal{O}_{X'}(n_2)\oplus \mathcal{O}_{X'}(-(n_1+n_2))$. This implies (i).
\\ (2)  If $\G\cong \mathcal{O}_{X'}(n)\oplus \M$ with $\M$ an indecomposable
 vector bundle of constant rank 2 then $\M$ must be absolutely indecomposable, since,
  if $l'/k'$ is a quadratic field extension, every line bundle over $X_{l'}$
 is already defined over $X'$. Moreover, we need
 ${\rm det}\,\M\cong\mathcal{O}_{X'}(-n)$. We get

 \[ \Pe_0\cong  \h om_{\mathcal{O}_{X'}}(\mathcal{O}_{X'}(n)\oplus \M,\mathcal{O}_{X'}\oplus\mathcal{O}_{X'}(m)\oplus\mathcal{O}_{X'}(-m))\cong\\
\left [\begin {array}{cc}
\mathcal{O}_{X'}(-n) & \M^\vee \\
\mathcal{O}_{X'}(m-n) &\mathcal{O}_{X'}(m)\otimes \M^\vee \\
\mathcal{O}_{X'}(-m-n) &\mathcal{O}_{X'}(-m)\otimes \M^\vee\\
\end {array}
\right ].
\]
Thus
\begin{align*}\Pe\otimes\mathcal{O}_{X'}\cong [\mathcal{O}_{X'}(-n) \oplus \M^\vee \oplus\mathcal{O}_{X'}(m-n)
 \oplus\mathcal{O}_{X'}(m)\otimes \M^\vee
\oplus \mathcal{O}_{X'}(-m-n) \oplus \mathcal{O}_{X'}(-m)\otimes \M^\vee ] \\
\oplus
[\mathcal{O}_{X'}(n) \oplus \M \oplus \mathcal{O}_{X'}(-m+n) \oplus \mathcal{O}_{X'}(-m)\otimes \M
\oplus \mathcal{O}_{X'}(m+n) \oplus \mathcal{O}_{X'}(m)\otimes \M].
 \end{align*}
If $\M$ is not already defined over $X$ (since $^\omega\M\cong\M^\vee$,
$ tr_{k'/k}( \M)\cong tr_{k'/k}( \M^\vee)$) then we get (ii) as a possible $\sts$-module structure of $\Pe$.
\\
If $\M$ is already defined over $X$, that is $\M=\M_0\otimes
\mathcal{O}_{X'}$ then we get
\begin{align*} \Pe\cong[ \sts(-n) \oplus \M_0^\vee \oplus\sts(m-n) \oplus\sts(m)\otimes \M_0^\vee
\oplus\sts(-m-n) \oplus\sts(-m)\otimes \M_0^\vee ]\oplus \\
[\sts(n) \oplus \M_0 \oplus \sts(-m+n) \oplus\sts(-m)\otimes \M_0
\oplus\sts(m+n) \oplus\sts(m)\otimes \M_0]
\end{align*}
as a possible $\sts$-module structure of $\Pe$, i.e., (iii).
 \\ (3) If $\G$ is absolutely indecomposable with trivial determinant (and there exist such $\sts$-module!), then
  \[ \Pe_0\cong  \h om_{\mathcal{O}_{X'}}( \G,\mathcal{O}_{X'}\oplus \mathcal{O}_{X'}(m)\oplus\mathcal{O}_{X'}(-m))\cong\\
    \G^\vee\oplus\mathcal{O}_{X'}(m)\otimes\G^\vee\oplus\mathcal{O}_{X'}(-m) \otimes\G^\vee \]
    and
\[\Pe\otimes\mathcal{O}_{X'}\cong [ \G^\vee\oplus\mathcal{O}_{X'}(m)\otimes\G^\vee\oplus\mathcal{O}_{X'}(-m) \otimes
\G^\vee ]\oplus
[ \G\oplus\mathcal{O}_{X'}(-m)\otimes\G\oplus\mathcal{O}_{X'}(m) \otimes\G ]. \]
If $\G$ is already defined over $X$, that is $\G=\G_0\otimes
\mathcal{O}_{X'}$ then we get (iv).
\\
If $\G$ is not already defined over $X$ then (since $^\omega\G\cong\G^\vee$,
$ tr_{k'/k}( \G)\cong tr_{k'/k}( \G^\vee)$) we get (v):

\[\Pe\cong tr_{k'/k}( \G)\oplus\sts(m)\otimes  tr_{k'/k}( \G)\oplus\sts(-m)\otimes  tr_{k'/k}( \G). \]

We have thus proved the assertion using [AEJ1] and the fact that every line bundle over $X'$ is already defined over $X$ here
which excludes the case that $\G$ is indecomposable, but not absolutely so
 (there is no vector bundle $\G$ of rank 3 which is indecomposable, but not absolutely so, since every
$\G\cong  tr_{l'/k'}(\mathcal{O}_{X_{l'}}(n))$, $l'/k'$ a cubic field extension, decomposes)
.
\end{proof}

Note that there are probably other vector bundles $\E$ over $X'$ of rank 3 carrying a nondegenerate hermitian form, so that the above
are only examples and do not cover all possible cases.

\smallskip
{\it Acknowledgements:}
The author would like to acknowledge the support of the ``Georg-Thieme-Ged\"{a}chtnisstiftung'' (Deutsche Forschungsgemeinschaft)
during her stay at the University of Trento,
 and thanks the Department of Mathematics at the University of Trento for its hospitality.

\end{document}